\documentclass[12pt,twoside]{article}
\long\def\symbolfootnote[#1]#2{\begingroup%
\def\thefootnote{\fnsymbol{footnote}}\footnote[#1]{#2}\endgroup}

\usepackage{amsfonts, amssymb, amsmath, amsthm}
\usepackage[colorlinks=true, pdfstartview=FitV, linkcolor=blue,
            citecolor=blue, urlcolor=blue]{hyperref}
\usepackage[usenames]{color}
\definecolor{Red}{rgb}{0.7,0,0.1}
\definecolor{Green}{rgb}{0,0.7,0}
\usepackage{accents}
\usepackage{comment}

\newcommand{\dom}{\mathcal M}
\renewcommand{\L}{L_2 (\mathfrak U, L^2(\dom))}
\newcommand{\Ri}{\mathcal R}

\newcommand{\Ra}{\bar{\mathcal R}}

\renewcommand {\u}{\tilde u}
\renewcommand {\v}{\tilde v}

\title
{Martingale and Pathwise Solutions to the Stochastic Zakharov-Kuznetsov  Equation with Multiplicative Noise}
\author{Nathan Glatt-Holtz$^{\sharp}$,  Roger Temam${^\flat}$,  Chuntian Wang${^\flat}$,
}

\date{}

\numberwithin{equation}{section}

\newtheorem{thm}{Theorem}[section]
\newtheorem{lem}{Lemma}[section]
\newtheorem{prop}{Proposition}[section]

\newtheorem{deff}{Definition}[section]
\newtheorem{rem}{Remark}[section]

\begin{document}


\maketitle

\vskip-4mm

\centerline{\footnotesize{\it $^{\sharp}$ Institute for Mathematics and Its Applications } }
\vskip-1mm
\centerline{\footnotesize{\it University of Minnesota }}
\centerline{\footnotesize{\it  Minneapolis, MN 55455  }}
\vskip-1mm
\centerline{\footnotesize{\it email: \url{negh@ima.umn.edu}}}
\vskip-1mm
\centerline{\footnotesize{\it $^{\sharp}$ {Department} of Mathematics, Virginia Polytechnic and State University} }
\vskip-1mm
\centerline{\footnotesize{\it  Blacksburg, VA 24061  }}
\vskip-1mm
\centerline{\footnotesize{\it email: \url{negh@vt.edu}}}

\vskip2mm
\centerline{\footnotesize{\it ${^\flat}$ Department of Mathematics and The Institute for Scientific Computing and Applied Mathematics}}
\vskip-1mm
\centerline{\footnotesize{\it Indiana University, Bloomington, IN 47405}}
\vskip-1mm
\centerline{\footnotesize{\it email: \url{temam@indiana.edu }}}
\vskip 0mm
\centerline{\footnotesize{\it \,\, email: \url{wang211@umail.iu.edu}}}

\begin{center}
\large
\date{\today}
\end{center}

\vskip4mm

\tableofcontents

\newpage

\begin{abstract}
We study in this article  the stochastic Zakharov-Kuznetsov equation driven by a multiplicative noise.
We establish, in space dimensions two and three the global existence of martingale solutions,  and  in space dimension two the global pathwise uniqueness and the existence of pathwise solutions.
New methods are employed in the passage  to the limit on a
special type of boundary conditions and in the verification of the pathwise uniqueness of martingale solutions with a lack of regularity, where both difficulties arise due to the partly hyperbolic feature of the model.
%
%
\end{abstract}
{\noindent \small
 {\it \bf Keywords: Zakharov-Kuznetsov equation,  Korteweg-de Vries equation,  stochastic fluid dynamics.\\}


\section{Introduction}
\label{sec:introduction}

We consider the stochastic Zakharov-Kuznetsov equation subject to 
multiplicative random noise
\begin{equation}\label{eq1111}
   d u+ (\Delta u_x  + c u_x + u u_x )\, d t=f\, dt + \sigma(u) \, d W(t),
  \end{equation}
where $u=u(x,  x^\perp, t)$,  $x^\perp =y$ or $x^\perp=(y,z)$ in a limited domain $ \{  (x, x^\perp),\, 0<x<1, \,  x^\perp\in (-\pi /2, \pi /2)^d, \, d=1, 2    \}$.

The deterministic ZK equation  describes the propagation of nonlinear ionic-sonic waves in a plasma submitted to a magnetic field directed along the $x$-axis.     It has been derived formally in a long wave, weakly nonlinear regime from  the Euler-Poisson system in \cite{ZK},
 \cite {LSp} and \cite{LLS} (see also \cite{bonaparitchard1} and \cite{bonapritchatd2} for more general physical backgrounds).
When $u$ depends only on $x$ and $t$, the  ZK equation reduces to the classical Korteweg-de Vries (KdV) equation.
Recently the ZK equation has caught considerable attention (see e.g. (\cite{Fam2}), \cite{Fam3},  \cite{LT},  \cite{DoroninLarkin}, \cite{BaykovaFaminskii},  \cite{SautTemam} and \cite{SautTemamChuntian}),  not only because it is closely related with the physical phenomena but also because it is the start to explore more general problems that are partly hyperbolic (such as the inviscid primitive equations).
%
%

To capture the phenomena similar to those of more realistic fluid systems, random waves   have attracted interests nowadays. For example, the stochastic KdV equation has been studied extensively (see
   \cite{bouardannearnaud}, \cite{gaowabao} and \cite{debusschepinterms})), where the main focus are on  Wick-type SPDEs (\cite{zhangzhang} and \cite{Liu}) and
 exact solutions of the stochastic KdV equation under additive noise (\cite{hermanrose}).
However, to the best of our knowledge, there has been no result so far for the  stochastic  ZK equation driven by multiplicative noise. In particular, although (\ref{eq1111}) bears some similarity to the model studied in e.g. \cite{DebusscheGlattHoltzTemam1}, the situation is in fact totally different and we can not use the results derived  there; specifically, the operator  $Au:= \Delta u_x + cu_x$ is neither symmetric nor coercive in our case (see Remark \ref{eq1111} for details).

In the present article, we extend the results of global existence (in space dimensions two and three) and uniqueness (in space dimension two) of weak solutions  in \cite{SautTemamChuntian} to the stochastic case.   This initial program of well posed-ness will serve as the foundation for the investigation of the qualitative and quantitative properties of solutions in both the deterministic and stochastic cases and facilitate the comparison of these two models.
Note that here we have different notions of solutions, namely, the martingale solutions and the pathwise solutions. In the former notion, the stochastic basis is not specified in the beginning and is viewed as part of the unknown, while in the latter case, the stochastic basis is fixed in advance as  part of the assumptions. To pass from the martingale to the pathwise solutions,
 we apply the extension of  the  Yamada-Watanabe theorem (\cite{WatanabeYamada1971a}) to the infinite dimension (see \cite{GyongyKrylov1}), that is, the pathwise solutions exist whenever there exists a pathwise unique martingale solution.


One of the main novelties of this paper is the treatment to the boundary conditions, which are more complicated than the usual Dirichlet or periodic ones. Firstly, it is not clear whether  all the boundary conditions are still preserved after the
application of the  Skorokhod  embedding theorem (Theorem 2.4 in \cite{ZabczykDaPrato1}) since the underlying stochastic basis has been changed.  To  solve this problem, a measurability result concerning Hilbert spaces is developed (Lemma \ref{lemcompact1}).
Secondly, we have extended the trace results in \cite{SautTemamChuntian} to the stochastic setting by establishing the trace properties of the linearized ZK equation   depending on the probabilistic parameters (Lemma \ref{lineartracet} and \ref{A.2t}). This method can be used to deal with the non-conventional boundary conditions  in other circumstances in  future.

A further novelty is contained in the proof of the pathwise uniqueness (Section \ref{lacalp}). Difficulties arise with the derivation of  the energy inequality for the difference of the solutions due to the lack of regularity. Moreover, the method in the deterministic case (see \cite{SautTemamChuntian}) can not be adapted to the stochastic case if we just use  the  stochastic version of   the Gronwall lemma established in \cite{GlattHoltzZiane2} (see also \cite{MR}),  as  issues would arise in passage to the limit on the terms involving stopping times. We overcome this difficulty  by establishing a variant of the stochastic Gronwall lemma (Lemma \ref{lemsg1}), where we find that in certain situations  we can weaken the hypotheses so as to avoid the stopping times. This idea of making use of Lemma \ref{lemsg1} to deal with the insufficient regularity of the solutions  to verify the pathwise uniqueness, we believe, is not only suited for this model  but can also be applied to other classes of equations.

%
%
%
%
%

The present article is organized as follows. 
In Section \ref{section10} we make precise the assumptions on the problem  (\ref{eq1111}) and the stochastic framework. In Section \ref{regularize}, we study the parabolic regularization of equation (\ref{eq1111}) in  dimensions two and three. We show by Galerkin truncation the existence and uniqueness of the global pahwise solution $u^\epsilon$ which is sufficiently regular for the subsequent calculations and more importantly, we establish  the uniform estimates independent of $\epsilon$, which we use in Section \ref{subsection12} to develop the compactness argument when $\epsilon>0$ varies. Then  with application of the    Skorokhod  embedding  theorem, which leads to strong convergence of some subsequence, we obtain  the global existence of martingale solutions to (\ref{eq1111})    in  dimensions two and three.  In Section \ref{lacalp} we prove in  dimension two the pathwise uniqueness of martingale solutions and by the Gy\"ongy-Krylov  Theorem deduce the global existence of pathwise solution.

Finally in the Appendix, we present the generalized trace results, the measurability result and   the adapted  stochastic Gronwall lemma, among the other existing results used in the article.


%


\section{Stochastic ZK equation}
\label{section10}

We consider the stochastic ZK equation with multiplicative noise
\begin{equation}\label{eq111}
   d u+ (\Delta u_x  + c u_x + u u_x )\, d t=f\, dt + \sigma(u) \, d W(t),
  \end{equation}
evolving in a rectangular or parallelepiped domain, namely,  in $\dom=(0,1)_{x}\times(-\pi /2, \pi /2)^d,$ with $ d=1$ or $2$. In the sequel, we will use the notations $I_{x}=(0,1)_{x},$ $ I_{y}=(- \pi /2,  \pi /2)_{y}$ and $I_{z}=(- \pi /2, \pi /2)_{z}$.  Here $\Delta u=u_{xx}+\Delta^\perp u$, $\Delta^\perp u=u_{yy}$ or $u_{yy}+u_{zz}$. We assume that
$f$ is a  deterministic  function, and the white noise driven stochastic term  $\sigma (u) \, dW (t)$ is in general state dependent.
As in  \cite{SautTemamChuntian}, we assume the boundary conditions on $x=0, 1$
 to be
\begin{equation}\label{eq3}
 u\big|_{x=0}=u\big|_{x=1}=u_x\big|_{x=1}=0.
\end{equation}
For the  boundary conditions in the $y$ and $z$ directions,  we can choose  either the Dirichlet boundary conditions
\begin{equation}\label{eq107}
  u=0 \mbox { at } y=\pm\frac{\pi}{2}\,\, \,\left(\mbox{and } z=\pm\frac{\pi}{2}\right),
\end{equation}
or the periodic boundary conditions
\begin{equation}\label{eq1-26p}
 u\big|^ {y= \frac {\pi}{2}}_{y=-\frac {\pi}{2}}=u_y\big|^ {y= \frac {\pi}{2}}_{y=-\frac {\pi}{2}}=0\,\,\,\left(\mbox{and }u\big|^{z= \frac {\pi}{2}} _{z=-\frac {\pi}{2}}=u_z\big|^{z= \frac {\pi}{2}} _{z=-\frac {\pi}{2}}=0\right).
\end{equation}
The initial condition reads:
\begin{equation}\label{eq33}
 u(0)=u_0.
\end{equation}
 For the simplicity of the presentation, we will mostly study the Dirichlet case  (\ref{eq111})-(\ref{eq107}) and   (\ref{eq33}). We will just make some remarks concerning the closely related space periodic case when (\ref{eq107}) is replaced by (\ref{eq1-26p}).

\subsection{Stochastic framework}
\label{stochastic}

In order to define the term $\sigma(u) \, d W(t)$ in (\ref{eq111}), we  recall some basic notions and notations of stochastic analysis from \cite{DebusscheGlattHoltzTemam1}. For further details and background, see e.g. \cite{PrevotRockner}, \cite{FlandoliGatarek1},   \cite{Flandoli1}, \cite{Bensoussan1} and  \cite{ZabczykDaPrato1}.

To begin with we fix a stochastic basis
\begin{equation}\label{stochasticbasis}
\mathcal S:= \left( \Omega, \mathcal F, \{\mathcal {F}_t \}_{t\geq 0}, \mathbb P, \{ W_k \}_{k \geq 1} \right),
\end{equation}
that is a filtered probability space with $\{ W_k \}_{k \geq 1}$ a sequence of independent standard one-dimensional Brownian motions relative to $\{\mathcal {F}_t \}_{t\geq 0}$. In order to avoid unnecessary complications below we may assume that $\mathcal F_t$ is complete and right continuous (see \cite{ZabczykDaPrato1}).

We fix a separable Hilbert space $\mathfrak U$ with an associated orthonormal basis $\{ e_k \}_{k\geq 1}$. We may formally define $W$ by taking
$W= \sum _{k=1} ^ {\infty}W_k e_k$. As such $W$ is said to be a `cylindrical Brownian motion' evolving over $\mathfrak U$.

We next recall some basic definitions and properties of spaces of Hilbert-Schmidt operators.
 For this purpose we suppose that $X$  is any  separable Hilbert space with the associated norm and inner product written as $|\cdot|_X$, $\left (\cdot, \cdot \right)_X$. We denote by
\begin{equation*}
L_2(\mathfrak U, \, X)= \left \{  R \in  L ( \mathfrak U, \, X):\,\,\, \sum_k |R e_k|^2 _X < \infty  \right \},
\end{equation*}
the space of Hilbert-Schmidt operators from $\mathfrak U$ to $X$. We know that the definition of $L_2(\mathfrak U, \, X)$ is independent of the choice of the orthonormal basis $\{ e_k \}_{k\geq 1}$ in $X$.
 By endowing this space with the inner product $\left <R, T  \right > _{L_2 (\mathfrak U, X)} = \sum _k \left <Re_k, T e_k \right > _ X$, we may consider $L_2 (\mathfrak U, X)$ as itself being a Hilbert space. Again this scalar product can be shown to be independent of the orthonormal basis $\{ e_k \}_{k\geq 1}$.

%
%
%

We also define the auxiliary space $\mathfrak U_0 \supset \mathfrak U$ via
\begin{equation*}
\mathfrak U_0: = \left  \{ v =  \sum _{k \geq 0} a_k e_k:\,\, \sum _k \frac{a^2_k}{ k ^2} < \infty  \right    \},
\end{equation*}
 endowed with the norm $|v|^2_{\mathfrak U_0}:= \sum _k a_k^2 / k^2$, $v = \sum _k a_k e_k$. Note that the embedding of $\mathfrak U \subset \mathfrak U_0$ is Hilbert-Schmidt.
  Moreover, using standard martingale arguments combined with the fact that each $W_k$ is almost surely continuous (see \cite{ZabczykDaPrato1}) we obtain that, for almost every $\omega \in \Omega$, $W(\omega)\in \mathcal C ([0, T], \mathfrak U_0)$.

Given an $X$-valued predictable process $ \Psi \in L^2(\Omega; \, L^2 ((0, T), \, L_2 (\mathfrak U , \, X))$, one may define the It\=o stochastic integral
\begin{equation*}
M_t:= \int^t_0 \Psi\, dW= \sum_k \int^t_0 \Psi_k \, dW_k,
\end{equation*}
as an element in $\mathcal M^2_X$, that is the space of all $X$-valued square integrable martingales. In the sequel we will use the Burkholder-Davis-Gundy inequality which takes the form
\begin{equation}\label{burkholder}
 \mathbb E \left( \sup_{0\leq s \leq T} \, \left |\int ^s_0 \Psi  d\, W(t)\right|_X^r \right)
\leq c_1 \, \mathbb E\, \left [ \left( \int ^T_0||\Psi  ||^2_{L_2 (\mathfrak U , \, X)}  \, d  t \right ) ^ {r/2} \right ],
\end{equation}
valid for any $r\geq 1$. Here $c_1$ is an absolute constant depending only on $r$.

%

\subsection{Conditions imposed on $\sigma$, $f$ and $u_0$.}
\label{conditions}
 Given any pair of Banach spaces $\mathcal X_1$ and $\mathcal X_2$, we denote by $Bnd_u(\mathcal X_1, \, \mathcal X_2 )$, the collections of all continuous mappings
\begin{equation}\label{psi}
\Psi: \mathcal X_1 \rightarrow \mathcal X_2,
\end{equation}
such that
\begin{equation}\label{bnd}
||\Psi(u)||_{\mathcal X_2} \leq c_B ( 1 + ||u||_{\mathcal X_1} ), \quad u \in \mathcal X_1,
\end{equation}
for some constant $c_B$.
In addition, if
\begin{equation}\label{unif}
||\Psi(u)  -  \Psi(v)  ||_{\mathcal X_2} \leq c_U   ||u - v||_{\mathcal X_1} , \quad \forall u, \, v\in \mathcal X_1,
\end{equation}
for some constant $c_U$, we say that $\Psi \in Lip_u ( \mathcal X_1, \, \mathcal X_2  )$. In the sequel we will consider time dependent families of such mappings $\Psi=\Psi(t)$ and require that (\ref{bnd}) and (\ref{unif}) hold for $a.e.\,\,t$ with the same constants $c_B$, $c_U$ for all $t$'s under consideration.

We shall assume throughout the work that
\begin{equation}\label{bndu}
\sigma: [0, \, \infty ) \times L^2(\dom) \rightarrow L_2 (\mathfrak U, \,L^2(\dom)).
\end{equation}
Here $\mathfrak U$  and $ L_2 (\mathfrak U, L^2(\dom))$ are as introduced above. Moreover we assume that for $a.e.$ $t$,
\begin{equation}\label{sigma}
\sigma (t)\in Bnd_u ( L^2(\dom), \, \L  )\cap Bnd_u( \Xi_1 , \, L_2 (\mathfrak U,\,  \Xi_1)  ),
\end{equation}
and
\begin{equation}\label{unifsigma}
\sigma (t) \in Lip_u ( L^2(\dom), \, \L  ),
\end{equation}
where
\begin{equation}\label{xi1}
\Xi_1 := \left \{ u \in H^2(\dom) \cap H^1_0(\dom),\,\,u _x  \big|_{x=1}=0 \right  \}.
\end{equation}
When proving pathwise uniqueness of  martingale solutions and the existence of pathwise solutions in Section \ref{lacalp}, we will additionally suppose that for a.e. $t$,
\begin{equation}\label{sigmapa}
\sigma (t)\in Lip_u (L^2 (\dom) , \,  L_2 (\mathfrak U,  \, \Xi_1)   ).
\end{equation}
Furthermore in the sequel $\sigma$ is a measurable function of $t$ and all the corresponding norms of $\sigma(t)$ are essentially (a.e.) bounded in time.


%

Finally we state the assumptions for the initial condition $u_0$ and for $f$.  On the one hand, in Section \ref{subsection12}, where we consider only the case of martingale solutions, since the stochastic basis is an unknown of the problem, we will only be able to specify  $u_0$ as an initial probability measure $\mu_{u_0}$ on the space $L^2 (\dom)$ such that
\begin{equation}\label{martingaleinitial}
\int _{ L^2(\dom)} | u |_{{L^2 (\dom)}}^ 6 \, d\, \mu_{u_0}(u) < \infty,
\end{equation}
and we assume that $f$ is deterministic,
\begin{equation}\label{f}
f= f(x, x^\perp, t)\in  L^6 (0, T; \, L^2 (\dom)).
\end{equation}

On the other hand, for pathwise uniqueness and the existence of   pathwise solutions in Section \ref{lacalp},
where the stochastic basis $\mathcal S$ is fixed in advance we assume that,  relative to this basis, $u_0$ is an $L^2 (\dom)$-valued random variable such that
\begin{equation}\label{u0}
u_0 \in L^{7} (\Omega; \, L^2(\dom)) \mbox{ and $u_0$ is $\mathcal F_0$ measurable},
\end{equation}
and $f$ is deterministic,
\begin{equation}\label{fp}
f=f(x, x^\perp, t)\in  \, L^7 (0, T; \, L^2 (\dom)).
\end{equation}

\section{Regularized stochastic  ZK equation}
\label{regularize}


 As indicated above we consider the Dirichlet case, i.e.  (\ref{eq111})-(\ref{eq107}) and   (\ref{eq33}).  The domain is $\dom = I_x \times (-\pi/2, \pi/2)^d$, in $\mathbb{R} ^{d+1}$ with $d=1$ or $2$. In order to study this system,
 we will use a parabolic regularization of equation (\ref{eq111}), as in \cite{SautTemam}. That is, for $\epsilon >0$ ``small'', we consider the  stochastic parabolic equation of the $4$-th order in space:
%
%
\begin{equation}\label{10-11b}\displaystyle
\begin{cases}
&d u^\epsilon+ \left [\Delta u^\epsilon _x+ c   u^\epsilon_x + u^\epsilon u^\epsilon _x  + \epsilon \left( \frac{ \partial^4 u^\epsilon }{\partial x^4}  +
\frac{ \partial^4 u^\epsilon }{\partial y^4}  + \frac{ \partial^4 u^\epsilon }{\partial z^4}   \right )\right]dt =f^\epsilon \,dt + \sigma(u^\epsilon) \, d W(t),\\
&u^\epsilon (0)=u^\epsilon_0,
\end{cases}
\end{equation}
 supplemented with the  boundary conditions (\ref{eq3}),   (\ref{eq107}) and the additional boundary conditions
\begin{equation}\label{10-4}
u^\epsilon _{yy}\big|_{ {y=\pm \frac {\pi}{2}}  }=u^\epsilon _{zz}\big|_{ {z=\pm \frac {\pi}{2}}  }=0,
\end{equation}
\begin{equation}\label{10-2}
u^\epsilon _{xx}\big|_{x=0}  = 0.
\end{equation}
%

In the case of martingale solutions in Section  \ref{subsection12}, observing that  the space $L^2 (\Omega; \, \Xi_1)$ $\cap$ $L^{22/3} (\Omega; \, L^2 (\dom))$ is dense in $ L^6( \Omega; \, L^2 (\dom))$, we can use e.g. the Fourier series to construct an  approximate family $\{u^\epsilon_0\}_{\epsilon >0}$ which is $\mathcal F_0$ measurable, such that, as $\epsilon \rightarrow 0$:
 \begin{equation}\label{uepsilon0reg}
u^\epsilon_0 \in L^2 (\Omega;\, \Xi_1) \cap L^{22/3} (\Omega;\, L^2 (\dom)),
\end{equation}
\begin{equation}\label{initialc}
u^\epsilon_0 \rightarrow u_0 \mbox{ in } L^6( \Omega;\,  L^2 (\dom)).
\end{equation}
Similarly there exists a family of deterministic functions $\{f ^ \epsilon\}_{\epsilon >0}$ such that as $\epsilon \rightarrow 0$:
\begin{equation}\label{fep}
f^\epsilon \in  L^{22/3} (0, T; \, L^2 (\dom)),
 \end{equation}
\begin{equation}\label{initialfc}
f^\epsilon \rightarrow f \mbox{ in } L^6 (0, T; \, L^2 (\dom)).
\end{equation}

In the case  of pathwise solutions in Section \ref{lacalp}, in the same way we can deduce the existence of  the  approximate families $\{u^\epsilon_0\}_{\epsilon >0}$ and $\{f ^ \epsilon\}_{\epsilon >0}$ satisfying (\ref{uepsilon0reg}) and (\ref{fep}) respectively, and such that as $\epsilon \rightarrow 0$:
\begin{equation}\label{initialcp}
u^\epsilon_0 \rightarrow u_0 \mbox{ in } L^7( \Omega;  \, L^2 (\dom)),
\end{equation}
\begin{equation}\label{initialfcp}
f^\epsilon \rightarrow f \mbox{ in } L^7 (0, T; \, L^2 (\dom)).
\end{equation}

For notational convenience, as in \cite{SautTemamChuntian}, we  recast (\ref{10-11b}) in the form
  \begin{equation}\label{eq10-111}
  \begin{cases}
&d u^\epsilon=(- Au^\epsilon - B(u^\epsilon) - \epsilon \, L u^\epsilon  +f^\epsilon )\,d t + \sigma(u^\epsilon ) \, d W(t),\\
&u^\epsilon (0)=u^\epsilon_0,
\end{cases}
\end{equation}
where
\begin{equation}\label{opera}
\begin{split}
A u = \Delta u_ x + cu_x,\quad\quad\quad \quad\quad \quad \,\,\,\, &\forall\,\, u\in D(A),\\
 B(u, v) = u v_x   \in H^{-1}(\dom),\,\,\,\,\,\,\, \,\,\,\quad & \forall \,\,u \in L^2 (\dom),\,\, v\in H^1(\dom),\\
\,\,\,\,\,\,Lu =  u_{xxxx}   + u_{yyyy}  +  u_{zzzz}, \quad \,\quad & \forall \,\,u\in H^4(\dom).
\end{split}
\end{equation}
with  $D(A)=\left \{u\in L^2 (\dom): \,\,Au\in L^2 (\dom), \,\, u=0 \,\,\mbox{on}\,\, \partial\dom, u_x=0 \,\,at\,\, x=1 \right\}$. Note that a trace theorem proven in \cite{SautTemamChuntian} shows that if $u\in L^2 (\dom)$ and $Au \in L^2 (\dom)$ then the traces of $u$ on $\partial \dom$ and of $u_x$ at $x=1$ make sense.

\begin{rem}\label{intro}
As mentioned in the Introduction, although we can rewrite  (\ref{eq111}) as
 \begin{equation}\label{eq10-111p}
d u + ( Au +  B(u)) \, dt =  f  \,d t + \sigma(u ) \, d W(t),
\end{equation}
 which is similar to  the equation studied in \cite{DebusscheGlattHoltzTemam1}, the models are actually different.
 Indeed, the operator $A$ does not satisfy the assumptions in \cite{DebusscheGlattHoltzTemam1}; for example, $A$ is not symmetric. In fact,  for the adjoint $ A^*$ and  its domain $D(A^*)$, we have
 \begin{equation}\label{227}
 \begin{split}
D(A^*)&=\lbrace \bar u\in L^2 (\dom): \,\,A\bar u\in L^2 (\dom), \,  \bar u=0 \,\,on\,\, \partial\dom,\, \bar u_x=0 \,\,at\,\, x=0 \rbrace,\\
A^*\bar u&=-(\Delta \bar u_x+c \bar u_x), \,\,\,\bar u\in D(A^*).
\end{split}
\end{equation}
  For more details see  Section 2.3.2 in \cite{SautTemamChuntian}.

\end{rem}


\subsection{Definition of solutions}

We first introduce the necessary operators and functional spaces.
We will denote by $(\cdot, \cdot)$ and $|\cdot|$ the inner product and the norm of $L^2(\dom)$.
The space
$\Xi_1 $  defined in (\ref{xi1}) is  endowed with the scalar  product and  norm  $[ \cdot, \cdot ]_2 $, $ [ \cdot ]_2$:
\begin{align*}
[u  , v ]_2 =  (u_{xx}, v_{xx}) + (u_{yy}, v_{yy}) + (u_{zz}, v_{zz}),
\end{align*}
\begin{equation}\label{[]}
[u]_2^2  =    \left| u_{xx} \right|^2 + \left| u_{yy} \right|^2 +\left| u_{zz} \right|^2  ,
\end{equation}
which make it a Hilbert space. Note that since $|\Delta u| + |u|$ is a norm on $H^1_0 \cap H^2$ equivalent to the $H^2$-norm,
$[\cdot]_2$ is  a norm on $\Xi_1 $ equivalent to the $H^2$-norm.
Thanks to the Riesz theorem, we can associate to the scalar product $[ \cdot, \cdot ]_2 $ the isomorphism $\mathcal L$ from $\Xi_1$ onto ${\Xi_1}^\prime$, where $\mathcal L$ denotes the abstract operator corresponding to the differential operator $L$. Then considering the Gelfand triple $\Xi_1 \subset H:=  L^2 (\dom)\subset{\Xi_1}^\prime $, we introduce $\mathcal L^{-1} (H)$ the domain of $\mathcal L$ in $H$, which is the space
 \begin{equation}\label{spaceda}
\Xi_2  = \left \{ u \in \Xi_1 \cap H^4 (\dom),\, u _{yy}\big|_{ {y=\pm \frac {\pi}{2}}  }=u _{zz}\big|_{ {z=\pm \frac {\pi}{2}}  } = u _{xx }  \big|_{x=0}=0  \right  \}.
\end{equation}
The operator $\mathcal L^{-1}$ is self adjoint and compact in $H$. It possesses an orthonormal set of eigenvectors which is complete in $H$, and which we  denote by $\{\phi_i  \}_{i \geq 1}$. Note that all the $\phi_i$ belong to $\Xi_2$ which is the domain of $\mathcal L$ in $H$. Hence we have
\begin{equation*}
(\mathcal L u, v) = [ u  , v]_2, \,\,u \in \Xi_2,\, v \in \Xi_1.
\end{equation*}

We now introduce the following definitions.
\begin{deff}\label{defweaksolu}
(Global martingale solutions for the regularized ZK equation) Fix an $\epsilon >0$. For the case of martingale solutions, we only specify the measure $\mu_{ u^\epsilon_0}$ to be the probability measure of $u^\epsilon_0$ on  $\Xi_1$ which  satisfies
\begin{equation}\label{martingaleinitialb}
\int _{ L^2 (\dom)} | u |^ {22/3} \, d\, \mu _{u^\epsilon_0}(u) < \infty,
\end{equation}
\begin{equation}\label{martingaleinitialbb}
\int _{ \Xi_1 } | u|^ 2 \, d\,\mu _{u^\epsilon_0} (u) < \infty,
\end{equation}
and $f^\epsilon$ and $\sigma $ satisfy  (\ref{fep}), (\ref{sigma}) and (\ref{unifsigma}).


 A pair $(\tilde {\mathcal {S}}, \tilde u^\epsilon)$ is a global martingale solution to the regularized stochastic ZK equation  (\ref{10-11b})-(\ref{10-2}),   (\ref{eq3}) and    (\ref{eq107})  (in the Dirichlet case), if $ {\tilde {\mathcal S}}= ( \tilde\Omega, \tilde{ \mathcal F}, \{\tilde {\mathcal {F}}_t \}_{t\geq 0},  \tilde {\mathbb P}, \{ \tilde{ W}^k \}_{k \geq 1} )$ is a stochastic basis, and $ \tilde u^\epsilon (\cdot): \tilde  \Omega \times [0, \infty) \rightarrow  \Xi_1 $ is an $ \{ \tilde {\mathcal {F}}_t \} $ adapted process such that:
\begin{equation}\label{uepsiw}
\tilde u^\epsilon \in  L^{22/3} (\tilde \Omega; \, L^\infty (0, T; \, L^2 (\dom))) \cap L^2 (\tilde \Omega;\, L^\infty ([0, T]; \, \Xi_1) \cap L^2 (0, T; \, \Xi_2)),
\end{equation}
and
\begin{equation}\label{weakcontb}
 \tilde  u^\epsilon  (\cdot, \omega) \in \mathcal C ([0, T]; \, L^2_w (\dom)))\,\,\tilde {\mathbb P} - a.s.,
\end{equation}
where $L^2_w (\dom)$ is $L^2 (\dom)$ equipped with the weak topology,
and  the law of  $ \tilde  u^\epsilon (0)$ is $\mu _ { u^\epsilon _0}$, defined as above, i.e. $\mu _ { u^\epsilon _0}(E) = {\tilde {\mathbb P}} ( \tilde  u^\epsilon (0) \in E)$, for all Borel subsets $E$ of $ \Xi_1$, and finally $  \tilde  u^\epsilon $ almost surely satisfies
\begin{equation}\label{defsolb}
\tilde u^\epsilon (t) +   \int^t _0    (A\tilde u ^\epsilon + B (\tilde u^\epsilon ) + \epsilon L \tilde u ^\epsilon- f^\epsilon ) \,ds    = \tilde u^\epsilon (0)+  \int^t_0  \sigma (\tilde u^\epsilon) \, d\tilde{W},
\end{equation}
 as an equation in $L^2 ( \dom)$ for every $0 \leq t\leq T$.


\end{deff}


\begin{deff}\label{defstrongsol}
(Global pathwise solutions for the regularized ZK equation; Uniqueness)

Let $ {\mathcal S}:= ( \Omega, { \mathcal F}, \{ {\mathcal {F}}_t \}_{t\geq 0},  {\mathbb P}, \{ { W}^k \}_{k \geq 1} )$  be a  fixed stochastic basis and assume that  $u^\epsilon _0$, $\sigma$ and $f^\epsilon$ satisfy (\ref{uepsilon0reg}),    (\ref{sigma}), (\ref{unifsigma}) and (\ref{fep}).

(i) For any fixed $\epsilon >0$,  a random process $u^\epsilon$ is a global pathwise solution to      (\ref{10-11b})-(\ref{10-2}),   (\ref{eq3}) and    (\ref{eq107})   if $u^\epsilon$ is an $  {\mathcal {F}}_t  $ adapted process in $L^2 (\dom)$ so that (relative to the fixed-given-basis $\mathcal S$) (\ref{uepsiw})-(\ref{defsolb}) hold.

\vskip 1 mm

(ii) Global  pathwise solutions of   (\ref{10-11b})-(\ref{10-2}),   (\ref{eq3}) and    (\ref{eq107})  are said to be  global (pathwise) unique if given any pair of pathwise solutions $u^\epsilon$, $ v^\epsilon$ which coincide at $t =0$ on a subset $\Omega _0 $ of $ \Omega$, $\Omega _0  = \{ u^\epsilon (0) = v^\epsilon (0)  \}$, then
\begin{equation}\label{defsrongsolll}
\mathbb P \{  \textbf 1 _{ \Omega _0  } (u^\epsilon (t)  = v^\epsilon (t))\}  =1, \,\,0\leq t \leq T .
\end{equation}

%
%

\end{deff}

%

%
%

In the sequel, we will prove that there exists a unique global pathwise solution $u^\epsilon$ to  (\ref{10-11b})-(\ref{10-2}),   (\ref{eq3}) and    (\ref{eq107}), which is sufficiently regular for the calculations in Section \ref{subsection12} to be fully legitimate without  any need of further regularization. The existence of such a solution is basically classical (see e.g. \cite{ZabczykDaPrato1}, \cite{Flandoli1},  \cite{FlandoliGatarek1} and   \cite{DebusscheGlattHoltzTemam1}) for a parabolic problem like this, but we will make partly explicit the construction of $u^\epsilon$ because we need to see how the estimates depend  or  not on $\epsilon$.

\subsection{Pathwise solutions in dimensions $2$ and $3$}
\label{lacalpe}

With the above definitions, we can  state the main result of  section  \ref{regularize}:
\begin{thm}\label{uepsilon}
When $d =1,2$, suppose  that, relative to a fixed given stochastic basis $\mathcal S$, $u^\epsilon _0$ satisfies (\ref{uepsilon0reg}), and that $f^\epsilon$ and $\sigma$ satisfy   (\ref{fep}), (\ref{sigma}) and (\ref{unifsigma}), with $\epsilon > 0$ fixed arbitrary. Then  there exists a unique global pathwise solution $ u^\epsilon$ which satisfies (\ref{10-11b})
and the  boundary conditions (\ref{eq3}),  (\ref{eq107}), (\ref{10-4}) and  (\ref{10-2}).

\end{thm}
To prove this theorem,  
we first use a  Galerkin scheme to  derive the estimates   indicating  a compactness argument based on fractional Sobolev spaces and tightness properties of the truncated sequence. Then
by the  Skorokhod embedding theorem  (see Theorem 2.4 in \cite{ZabczykDaPrato1}, also \cite{Billingsley2} and \cite{Jakubowski})  we can pass to the limit in the Galerkin truncation and hence obtain
the global existence of   martingale solutions. Finally we deduce  the existence of  global pathwise solutions using
 pathwise uniqueness of martingale solutions  and  the Gy\"ongy-Krylov  Theorem (Theorem \ref{krilov} of the Appendix).
 Here we will only present in details the derivation of the estimates, which will be utilized in the subsequent investigations of the stochastic ZK equation in Section \ref{secunif}.

We start the proof of Theorem \ref{uepsilon} by introducing the Galerkin system.
 We define $P^n $ as the orthogonal  projector from $L^2 (\dom) $ onto $H^n$, the space  spanned by the first $n$ eigenfunctions  of $\mathcal L$, $\phi_1$, ..., $\phi_n$. We consider the Galerkin system as follows
\begin{equation}\label{10-11g}\displaystyle
\begin{cases}
&d u^ n + \left (A^n u^n  + B^n (u^ n ) \right ) \, dt  + \epsilon Lu^n \, dt =f ^ n\,dt + \sigma^n(u^n )  \, d W(t),\\
&u^n (0)=P_n u_0,
\end{cases}
\end{equation}
 where $u^n$ maps $\Omega \times [0, \, T]$ into $H^n$,  $A^n u^ n := P^n A u^ n$, $B^n (u^ n): = P^n B (u^ n)$, and $\sigma ^ n(u^ n ) := P^n (\sigma(u^ n) )$.  In  (\ref{10-11g}), $\epsilon$ being fixed,    we write for simplicity $u^n$ for $u^{\epsilon, n}$  and $f^n$ for $f^{\epsilon, n}$.

Equation (\ref{10-11g}) is equivalent to a system of $n$ stochastic differential equations for the components of $u^n$ and it is a classical result that there exists a unique regular pathwise  solution $u^n = u^{\epsilon, n}$ such that
\begin{equation}\label{unreg}
u^n  \in L^2(\Omega; \, \mathcal C  (0, T;\, H^n)).
\end{equation}



\subsubsection{Estimates independent of $\epsilon$ and $n$}
\label{uniformestin}
We first derive the following estimates  on $u^{\epsilon, n}$ independent of $\epsilon$ and $n$.
\begin{lem}\label{uniformestb}
With the same assumptions as in Theorem \ref{uepsilon},  if  $u_0^\epsilon$ and $f^\epsilon$ satisfy (\ref{initialc}) and  (\ref{initialfc}) respectively, then the following estimates hold for $u^n = u^{\epsilon, n}$ independently of $\epsilon$ and $n$:
\begin{equation}\label{uni2}
u_x^{\epsilon, n} \big|_{x=0} \mbox{ remains bounded in } L^2 (\Omega; \, L^2 (0, T; \, L^2 (I_{x^\perp}))),
\end{equation}
\begin{equation}\label{uni3}
\sqrt \epsilon u^{\epsilon, n} \mbox{ remains bounded in } L^2 (\Omega; \, L^2 (0, T; \, \Xi_1)),
\end{equation}
\begin{equation}\label{uni1}
u^{\epsilon, n}\mbox{ remains bounded in } L^{6} (\Omega; \, L^\infty (0, T; \, L^2 (\dom))).
\end{equation}

If we further assume that (\ref{initialcp}) and (\ref{initialfcp}) hold, then
\begin{equation}\label{uni11}
u^{\epsilon, n}\mbox{ remains bounded in } L^{7} (\Omega; \, L^\infty (0, T; \, L^2 (\dom))),
\end{equation}
 with the bounds in  (\ref{uni2})-(\ref{uni11}) independent of both $\epsilon$ and $n$.
\end{lem}

\noindent\textbf{Proof.} We start by applying the It\=o formula to     (\ref{10-11g}). This yields
 \begin{equation}\label{ito22b}
 d|u^n|^2 = 2   \,\left(u^n , \, \mathcal N^n( u^n) \right)\,dt + 2\left (u^n , \, \sigma ^ n(u^n)\,dW(t)\right )  +  || \sigma ^ n(u^n ) || ^2 _{\L  }\, dt,
 \end{equation}
 where   $\mathcal N^n( u^n):= - A^nu^n - B^n(u^n) - \epsilon \, L u^n +f^n$, and  
  \begin{equation}\label{nn}
   \begin{split}
 (u^n , \, \mathcal N^n( u^n) ) &= - (u^n,  A^nu^n )  -(u^n, B^n(u^n)) - \epsilon (u^n,\, L u^n) +(u^n, f^n)\\
& = - (u^n,  Au^n )  -(u^n, B(u^n)) - \epsilon (u^n,\, L u^n) +(u^n, f^n).
\end{split}
 \end{equation}
To compute the right-hand side of (\ref{nn}) we remember that by (\ref{unreg}), for every $t$, $u^n (t) \in H^n= span\,\, (\phi_1, ..., \phi_n) \in \Xi_2$ a.s.,  since all the $\phi_i$ belong to $\Xi_2$. Hence in particular $u^n(t)$ satisfies the boundary conditions in  (\ref{spaceda}).
We drop the super index $n$ for the moment and  perform the following calculations a.s. exactly as in \cite{SautTemamChuntian}:
 \begin{align*}
   & \bullet\,\, (Au, \,u ) = \int _{\dom} \Delta u_x \, u \, d \dom +  c \int _{\dom}u_x  u \, d \dom \\
  &\quad \quad\quad \quad \quad = \int_{\partial \dom}\dfrac{ \partial u_x}{\partial n } u \, d \,\partial \dom  -\int_{\dom} \nabla u_x \cdot \nabla u \, d \dom +c \int_{\dom} \dfrac{ \partial}{\partial x }\left (  \dfrac{  u^2}{2 } \right) \, d\dom\\
   &\quad \quad\quad \quad \quad = - \int_{\dom} \dfrac{ \partial}{\partial x }\left (  \dfrac{  | \nabla u|^2}{2 } \right) \, d \dom+ \dfrac{c}{2} \int_{I_{x^\perp}}  u^2 \big|_{x=0}^{x=1}  \, d I_{x^\perp}\\
   &\quad \quad\quad \quad \quad = -\dfrac{1}{2} \int _{I_{x^\perp}} | \nabla u|^2 \big| ^ {x=1} _{x=0} \, d I_{x^\perp} =
   \frac {1} {2} \big| u_x \big|_{x=0} \big|_{L^2\left(I_{x^\perp}\right)}^2, \\
  & \bullet\,\, (  B(u), \, u)=0,\\ 
         & \bullet\,\,    \epsilon ( Lu,  \, u) =\epsilon \int_{\dom}  \left ( \left| u_{xx}\right |^2 + \left| u_{yy}\right|^2 +
\left|u_{zz}\right|^2 \right) \,d \dom.
\end{align*}
Hence by (\ref{ito22b}) we find
\begin{equation}\label{itointb}
\begin{split}
d\, | u^n |^2 + & \left ( \left|u^n_x \big|_{x=0} \right|^2_{L^2\left(I_{x^\perp}\right)} + 2 \epsilon [u^n]^2_2 \right) d t  \\
&\,\,\,\,\,\,\,\,\,\,=\, 2   (f^n,\,  u^n)\, dt +  || \sigma ^ n(u^n) || ^2 _{\L } \, dt+2 \left( u^n, \,\sigma ^ n(u^n)\,dW(t)\right).
\end{split}
\end{equation}
Integrating both sides from $0$ to $s$ with $0 \leq s \leq r \leq T$, taking the supremum over $[0, r]$, we have
\begin{equation}\label{eq16b1}
\begin{split}
\sup_{0\leq s \leq r} \, |u^n (s) |^2 + &  \int^r_0  \left ( \left| u^n_x \big|_{x=0} \right|^2_{L^2\left(I_{x^\perp}\right)} + 2 \epsilon [u^n]^2_2 \right )\, dt\\
\leq \, &  |u^n_0|^2  +  2\,  \int^r_0   \left|(f^n,\,  u^n)\right|\, dt +  \int^r_0 || \sigma^n (u^n) || ^2 _{\L } \, dt \\
&+ 2 \sup_{0\leq s \leq r} \,\left | \int ^s_0 \left (u^n, \, \sigma^n  (u^n) \,dW(t)\right ) \right | .\end{split}
\end{equation}
Raising both sides to the power $p/2$ for $p\geq 2$, then taking expectations,
   we obtain with the Minkowski inequality and Fubini's Theorem
%
\begin{equation}\label{eq16b}
\begin{split}
 \mathbb  E \sup_{0\leq s \leq r} \, |u^n (s) |^p
\lesssim &\, \mathbb  E  |u^n_0|^p  +  2\,\mathbb E \, \int^r_0   \left|(f^n,\,  u^n)\right|^{p/2}\, dt \\
&\,+  \mathbb  E \int^r_0 || \sigma^n (u^n) || ^p _{\L } \, dt \\
&\,+ 2\mathbb  E \left (\sup_{0\leq s \leq r} \left | \int ^s_0 \left (u^n, \, \sigma^n  (u^n) \,dW(t)\right ) \right | \right)^{p/2},
\end{split}
\end{equation}
where $\lesssim$ means $\leq $ up to an $absolute$ multiplicative constant. Here and below  $c^\prime$ indicates  an absolute constant, whereas $\eta$,  $\kappa$, and the $\kappa_i$ indicate constants depending on the data $u_0$, $f$, etc. These constants may be different at each occurence.
We estimate the terms  on the right-hand side of (\ref{eq16b}) a.s. and for a.e. $t$:
\begin{align*}
& \left|(f^n,\,  u^n)\right|^{p/2}  \leq |f^n|^{p/2} | u^n| ^{p/2}  \leq     |u^n|^p + | f^n|^p,\\
&  || \sigma^n (u^n) || ^p _{\L } \leq  \mbox{ by (\ref{sigma}) }    \leq c_B^p (|u^n |+ c^\prime )^{p} \lesssim |u^n |^p+ c^\prime  ;
\end{align*}
%
 for the stochastic term, we use the Burkholder-Davis-Gundy inequality (see (\ref{burkholder})) and (\ref{sigma})
  \begin{equation*}\label{bdg}
  \begin{split}
 \quad \quad & \mathbb  E   \left (\sup_{0\leq s \leq r} \left | \int ^s_0 \left (u^n, \, \sigma^n  (u^n) \,dW(t)\right ) \right | \right)^{p/2} \\
&\quad\quad\quad\quad\quad\leq \mathbb E \sup_{0\leq s \leq r} \left | \int ^s_0 \left (u^n, \, \sigma^n  (u^n) \,dW(t)\right ) \right | ^{p/2}  \\
&\quad\quad\quad\quad\quad \leq c _1\, \mathbb E\, \left [ \left (\int ^r_0 |u^n |^{2} \, ||\sigma^n (u^n)||^2 _{\L }  \, d  t \right ) ^ {p/4} \right ]\\
&\quad\quad\quad\quad\quad  \lesssim \mathbb  E \left [ \left (  \sup_{0\leq s \leq r}|u^n |^{2}   \int ^r_0 1 + |u^n |^{2}\, d  t \right ) ^ {p/4} \right ] \\
&\quad\quad\quad\quad\quad  \leq \dfrac{1}{2} \, \mathbb  E   \sup_{0\leq s \leq r}|u^n |^{p} + c^{\prime }   \,  \mathbb  E  \int ^r_0 |u^n |^{p}\, d  t +c^\prime.
\end{split}
\end{equation*}
%
%
Applying the above estimates to (\ref{eq16b}), we obtain
%
 \begin{equation}\label{applyingb}
\dfrac{1}{2}\,  \mathbb E\sup_{0\leq s \leq r} \, |u^n(s) |^p  \leq \mathbb   E \,  |u^n _0|^p +   c^\prime \, \mathbb  E\int_0^r  |u^n(t)|^p \, d t+ \mathbb   E \int ^r_0  |f^n(t)|^p\,d t + c^\prime.\footnote{Note that here $f^{\epsilon, n}$ is actually independent of $\omega\in \Omega$ and the symbol $\mathbb E$ in front of the corresponding term is not needed. However in Section \ref{pathwisesolutions}   we will use another version of this calculation in which $f^{\epsilon, n}$ is replaced by $g^\epsilon$ which depends on $\omega$; hence we leave $\mathbb E$ in front of the term involving $f^{\epsilon, n}$ in view of the calculations in Section \ref{pathwisesolutions}.}
\end{equation}
Since $\mathbb  E\int_0^r  |u^n(t)|^p \, d t \leq  \int_0^r \,  \mathbb E \, \sup_{0\leq l\leq t} \, |u^n (l) |^p \, d t$, setting $ \mathbb E\, \sup_{0\leq s \leq r} \, |u^n(s) |^p =: U(r)$,  with (\ref{applyingb}) we deduce
\begin{equation*}
 U(r) \leq   U(0)  + c^\prime \int_0^r \, U(t) d t+ \, \int^r_0\, \mathbb E\,|f^n(t)|^p\,d t + c^\prime,
\end{equation*}
for every  $0\leq r\leq T$.
Hence  applying  the (deterministic) Gronwall lemma, we obtain for $p\geq 2$,
\begin{equation}\label{eq21bb}
\mathbb  E\sup_{0\leq r \leq T} \, |u^n (r) |^{p}
\lesssim  \mathbb E  |u^n _0|^p +    \mathbb   E \int ^T_0  |f^n(t)|^p\,d t + c^\prime .
\end{equation}
Letting $p=6$, thanks to (\ref{initialc}) and (\ref{initialfc}), we deduce that
\begin{equation}\label{eq21b}
\mathbb  E\sup_{0\leq r \leq T} \, |u^n (r) |^{6}
\leq \kappa_1,
\end{equation}
for a constant $\kappa_1$ depending only on $u_0$, $f$, $T$ and $\sigma$, and independent of $\epsilon$ and $n$; this  implies (\ref{uni1}).  Similarly, setting $p=7$ in (\ref{eq21bb}),  we infer (\ref{uni11}) from (\ref{initialcp}) and (\ref{initialfcp}). Finally, setting $p=2$ in (\ref{applyingb}), along with (\ref{eq16b1}) we obtain (\ref{uni2}) and (\ref{uni3}). \qed


\subsubsection{Estimates dependent on $\epsilon$}

We now derive estimates independent of $n$ only, that is, valid for fixed $\epsilon$.
\begin{lem}\label{uniformestb2}
With the same assumptions as in Theorem \ref{uepsilon},  the following  estimates hold for $u^n = u^{\epsilon, n}$,  as  $n \rightarrow \infty$ and $\epsilon>0$ remains fixed:
\begin{equation}\label{uni0b}
 u^{\epsilon, n}\mbox{ remains bounded in } L^{22/3} (\Omega; \, L^\infty (0, T; \, L^2 (\dom))),
\end{equation}
\begin{equation}\label{uni1b}
u^{\epsilon, n}\mbox{ remains bounded in } L^2 (\Omega; \, L^\infty  (0, T; \,\Xi_1)),
\end{equation}
\begin{equation}\label{uni4b}
u^{\epsilon, n} \mbox{ remains bounded in } L^2 (\Omega; \, L^2 (0, T; \, \Xi_2)).
\end{equation}
\end{lem}
\noindent\textbf{Proof.}
Setting $p=22/3$ in (\ref{eq21bb}), we infer (\ref{uni0b}) from (\ref{uepsilon0reg}) and (\ref{fep}).

Returning to (\ref{10-11g}), we apply  the  It\=o formula to (\ref{10-11g}) and  obtain an evolution equation for the $\Xi_1$ norm:
 \begin{equation}\label{eqdubb}
 d[u^n]_2^2 = 2   \,\left(L u^n , \, \mathcal N( u^n)\right)\,dt + 2\left ( L u^n , \, \sigma ^ n(u^n)\,dW(t)\right)  +  || \sigma ^ n(u^n ) || ^2 _{L_2 (\mathfrak U,\, \Xi_1)   }\, dt,
 \end{equation}
where $\mathcal N^n( u^n)$ has been  defined before.
Similar to (\ref{nn}) we have  a.s. and for a.e. $t$:
\begin{align}\label{un}
\left (L u^n , \, \mathcal N^n( u^n)\right )
 & = - (Lu^n,  Au^n )  - \left (Lu^n, B(u^n) \right) - \epsilon |L u^n|^2 +(Lu^n, f^n).
 \end{align}
 By (\ref{eqdubb})
 and (\ref{un})
 we deduce
 \begin{equation}\label{eqdubbb}
 \begin{split}
 d[u^n]_2^2 + \epsilon |L u^n|^2 \,dt  = & - 2 (Lu^n,  Au^n )  \,dt  -2 (Lu^n, B(u^n)) \,dt  + 2(Lu^n, f^n)\,dt\\
&+ 2\left ( L u^n , \, \sigma ^ n(u^n)\,dW(t) \right)  +  || \sigma ^ n(u^n ) || ^2 _{L_2 (\mathfrak U, \, \Xi_1)   }\, dt.
\end{split}
 \end{equation}
 Integrating both sides from $0$ to $s$ with $0 \leq s  \leq T$, taking the supremum over $[0, T]$,  then taking expectations, we arrive at
\begin{equation}\label{eq16bbb}
\begin{split}
\mathbb E\sup_{0\leq s \leq T} \, [u^n]_2^2 + \epsilon\, \mathbb  E \int^T_0 |L u^n|^2 \,dt
\leq&\,   \mathbb E\, |u^n_0|^2  +  2 \, \mathbb E \int^T_0 \left| (Lu^n,  Au^n ) \right| \,dt \\
&\,\,\, +  2 \, \mathbb E \int^T_0 \left| (Lu^n, B(u^n))\right| \,dt +  2 \, \mathbb E \int^T_0 \left|(Lu^n, f^n) \right| \,dt \\
&\,\,\,+2 \, \mathbb E \sup_{0\leq s \leq T} \, \int ^s_0 \left ( L u^n , \, \sigma ^ n(u^n)\,dW(t) \right)\\
& \,\,\,+ \mathbb  E \int^T_0 || \sigma ^ n(u^n ) || ^2 _{L_2 (\mathfrak U, \,\Xi_1)   }\, dt.
\end{split}
\end{equation}
We estimate a.s. each term on the right-hand side of (\ref{eq16bbb});  we emphasize that the estimates depend on $\epsilon$ but not on $n$ or on  $\omega \in \Omega$:
\begin{align*}
\left |(Lu^n,  Au^n ) \right| &\leq |Lu^n |   | u^n|_{H^3 (\dom)}
\leq |Lu^n | |u^n|^{1/4}_{H^2}|u^n|^{3/4}_{H^4}
\leq |Lu^n | ^{7/4}  [u^n]_2^{1/4}
\leq \frac {\epsilon }{4}  |Lu^n | ^2 + \eta({\epsilon}) [u^n]_2^2,
\end{align*}
where $\eta({\epsilon})$ depends on $\epsilon$.
For the term  $\left |(Lu^n, B(u^n)) \right |$, we first estimate a.s. $| B(u^n)|$ in  dimension three:
\begin{equation}\label{buuq}
\begin{split}
\left |u^n u^n _x\right| &\leq |u^n|_{H_0^1 (\dom) }^{3/2}|u^n|^{1/2}_{H^2(\dom)}\\
&\leq (\mbox{by interpolation in dimension three, }|u^n|_{H_0^1 (\dom) } \lesssim |u^n| ^{3/4} |u ^n|^{1/4}_{H^4 (\dom)},\\
 & \,\,\,\,\,\,\mbox{\,\,\,\,\,\, and } |u^n|_{H^2(\dom) }\lesssim |u^n| ^{1/2} |u ^n|^{1/2}_{H^4 (\dom)}) \\
&\lesssim |u^n |^{9/8} |u^n |^{3/8}_{H^4 (\dom)} |u^n |^{1/4} |u^n |^{1/4}_{H^4 (\dom)}\\
&\lesssim |u^n|^{11/8} |L u^n |^{5/8}.\quad\quad\quad\quad\quad\quad\quad\quad\quad\quad\quad\quad\quad\quad\quad\quad\quad\quad\quad
\end{split}
\end{equation}
Hence
\begin{align*}
\left |(Lu^n, B(u^n))\right | &  \leq | B(u^n)|\,|Lu^n|
\leq (\mbox{by } ( \ref{buuq}))
\lesssim |u^n|^{11/8} |L u^n |^{13/8}
\lesssim  \frac{\epsilon}{4} |L u ^n|^{2} + \eta (\epsilon)|u^n|^{22/3} ,
\end{align*}where $\eta({\epsilon})$ depends on $\epsilon$.
%
For the stochastic term, we have
\begin{align*}
 \mathbb E \sup_{0\leq s \leq T} \left | \int ^s_0 ( L u^n , \, \sigma ^ n(u^n)\,dW(t)) \right|
&\leq  (\mbox{by the  Burkholder-Davis-Gundy inequality (\ref{burkholder})})\\
& \leq c_1 \, \mathbb E\, \left [ \left (\int ^T_0 | L u^n|^{2} \, ||\sigma^n (u^n)||^2 _{\L }  \, d  t \right ) ^ {1/2} \right ]\\
&  \leq c_1 c_B^2 \, \mathbb E\, \left [ \left (\int ^T_0 | L u^n|^{2} \,  (  1 + |u^n|^2  ) \, d  t \right ) ^ {1/2} \right ]\\
&     \leq \eta(\epsilon)\, \mathbb  E   \sup_{0\leq s \leq T}|u^n |^{2} + \frac{\epsilon}{4} \, \mathbb  E  \int ^T_0 | L u^n|^{2}\, d  t+ \eta(\epsilon),
\end{align*}where $\eta(\epsilon)$ depends on $\epsilon$.

For the term $\mathbb  E \int^T_0 || \sigma ^ n(u^n ) || ^2 _{L_2 (\mathfrak U,\,\Xi_1)   }\, dt$, we infer from (\ref{sigma}) that
\begin{equation*}
\mathbb  E \int^T_0 || \sigma ^ n(u^n ) || ^2 _{L_2 (\mathfrak U, \,\Xi_1)   }\, dt \lesssim \mathbb  E\int_{0}^{T}  \left (1+ [u^n]_2^2 \right ) \, d t.
\end{equation*}

Collecting all the above estimates, along with (\ref{eq16bbb}) we deduce
\begin{equation}\label{final}
\begin{split}
&\dfrac{1}{2}\,  \mathbb E\sup_{{0}\leq s \leq{T}} \, [u^n (s)]_2^2 + \frac{\epsilon}{4} \, \mathbb  E \int^{T}_{0}  |Lu^n|^2 \, dt\\
&\quad\quad\quad\quad\quad\quad \leq \mathbb   E \, [u^n_0]_2^2 + \eta(\epsilon)\, \mathbb  E\int_{0}^{T}  [u^n(t)]_2^2 \, d t +   \eta (\epsilon) \, \mathbb  E\int_{0}^{T}    |u^n(t)|^{22/3} \, d t\\
&\quad\quad\quad\quad\quad\quad\,\,\,\,\,\,\,\,\,\,\,+\eta(\epsilon) \, \mathbb  E   \sup_{0\leq s \leq T}|u^n |^{2} + \eta(\epsilon)\,  \mathbb   E \int_{0}^{T}  |f^n(t)|^2\,d t + \eta(\epsilon) .
\end{split}\end{equation}
Hence we can apply (\ref{uni0b}) and (\ref{uni3}) to (\ref{final}),   and we  obtain (\ref{uni1b}) and (\ref{uni4b}).  Thus we have completed the proof of Lemma \ref{uniformestb2}. \qed

\subsubsection{Estimates  in fractional Sobolev spaces.}

We will  apply the
 compactness result based on fractional Sobolev spaces in Lemma \ref{lemcompact1} (of the Appendix)   with
\begin{equation}\label{y_1}
\mathcal Y :=   L^2 (0, T; \, H^1_0 (\dom)) \cap  W^{\alpha, 2} (0, T;\,\Xi_2^\prime),\quad 0<\alpha < \dfrac {1}{2},
\end{equation}
where $\Xi_2^\prime$ is the dual of $\Xi_2$ relative to $L^2 (\dom)$. For that purpose  we will  need  the  following   estimates on fractional derivatives of $u^{\epsilon, n}$.

%
\begin{lem}\label{uniformestib}
With the same assumptions as in Theorem \ref{uepsilon}, we have
\begin{equation}\label{eub}
 \mathbb E |u ^{\epsilon, n} |_{\mathcal Y }  \leq  \kappa_2 (\epsilon),
\end{equation}
\begin{equation}\label{nostochastictermb}
\mathbb E \left |   u^{\epsilon, n} (t)  - \int^t _0  \sigma^n (u^{\epsilon, n}  ) \, d W(s) \right | ^2 _  {H^1  (0, T;\, \Xi_2^\prime   )}   \leq \kappa _3 ,
\end{equation}
\begin{equation}\label{nostochastictermmm}
\quad\,\,\,\,\,\,\,\,\mathbb E \left |\int^t _0  \sigma(u^{\epsilon, n}  ) \, d W(s)  \right |^2_{W^{\alpha, 6} (0, T;\, L^2 (\dom ))} \leq  \kappa _4, \,\,\,\,\,\,\forall \,\, \alpha < \dfrac {1}{2},
\end{equation}
where  $\kappa_2 (\epsilon)$ is independent of $n$  (but  may depend on $\epsilon$ and other data), while $\kappa_3$ and $\kappa_4$ depend only on $u_0$, $f$, $T$ and $\sigma$, and are independent of $\epsilon$ and $n$.
\end{lem}

\noindent \textbf{Proof.} (\ref{10-11g}) can be written as
 \begin{equation}\label{eq10-111intb}
  \begin{split}
 u^n (t)=& u^n_0   - \int^t _0 A^n u^n  ds  -\int^t _0 B^n(u^n) ds \\
 &  -\epsilon \int^t _0  L u ^n  ds +\int^t _0f^n ds  + \int^t _0   \sigma ^n(u^n )  d W(s)\\
:=& J^n_1 +  J^n_2 +  J^n_3 +  J^n_4 +  J^n_5 + J^n_6 .\end{split}
\end{equation}

For  $J^n_2$, 
fixing  $u^\sharp\in D(A^\ast) $ we have a.s. and for a.e. $t$
 \begin{align*}\label{a.s.}
 \left |(A^nu^n, \, u^\sharp)\right | =   \left |(u^n, \, A^\ast P^n u^\sharp) \right |    \leq   |u^n| \left | P^n u^\sharp \right |_{D(A^\ast)} \leq (\mbox{since $ \Xi_2 \subseteq D(A^\ast)$}) \leq  |u^n| |   u^\sharp|_{\Xi_2}   .
%
\end{align*}
Hence
    \begin{equation}\label{au1}
   |A^nu^n |_{\Xi_2^\prime}  \lesssim |u^n|.
   \end{equation}
%
With (\ref{au1}) and (\ref{uni1}) we obtain
%
\begin{equation}\label{j2}
\mathbb  E\, |J^n _2|^{{6}}_{ W^{1, 6} (0, T; \, \Xi_2^\prime)} \mbox{ is bounded independently of  $n$ and $\epsilon$}.
\end{equation}

For  $J^n_3$,
%
 firstly  we observe that $\forall \,\,u^\sharp\in \Xi_2 $ (dropping the  super index $n$ for the moment),
\begin{equation}\label{linf}
\begin{split}
 \left | (B(u),  \,  u^\sharp) \right | & = \left |  \int_{\dom}\dfrac{\partial }{\partial x } ( \dfrac {u ^{2}}{ 2 } ) \,u^\sharp \,d\,\dom \right | =   \dfrac{1}{2} \left | \int_{\dom}u ^{2} u^\sharp_x  \,d\,\dom \right |\\
& \leq  \dfrac{1}{2} |u|^2 |u^\sharp_x |_{L^\infty (\dom)}\\
& \leq  (\mbox{with } H^3 (\dom)\subset L^\infty (\dom) \mbox{ in  dimension }3 )\\
&\lesssim  |u|^2  |u^\sharp_x |_{H^3  (\dom)}\\
&\lesssim  |u|^2  |u^\sharp |_{\Xi_2 };
\end{split}
\end{equation}
 hence
\begin{equation}\label{bu1}|(B^n (u^n), \, u^\sharp)|
 = |(B (u^n), \, P^n u^\sharp)|
 \lesssim  |u^n|^2  |P^n u^\sharp |_{\Xi_2 }
 \leq |u^n|^2  | u^\sharp |_{\Xi_2 },
\end{equation}
which implies that $|B^n(u^n) |_{ \Xi_2^\prime}\lesssim  |u^n|^2$.
%
This along with  (\ref{uni1})   implies that
\begin{equation}\label{j33}
\mathbb  E |B^n( u^n) |^2_{ L^2 (0, T; \, \Xi_2^\prime)}  \mbox{ is bounded independently of $n$ and $\epsilon$},
\end{equation}
and hence
\begin{equation}\label{j3}
\mathbb  E |J ^n _3 |^2_{ H^1(0, T; \, \Xi_2^\prime)}  \mbox{ is bounded independently of  $n$ and $\epsilon$}.
\end{equation}

For  $J^n_4$, we have, $\forall \,\,u^\sharp \in \Xi_2$, $\left |\left (Lu^n, \, u^\sharp\right ) \right | =\left |  \left (u^n, \, L u^\sharp\right ) \right | \leq | u^n| |L u^\sharp |$.
%
Hence $|Lu^n| _{\Xi_2^\prime} \lesssim | u^n|$. Thus
\begin{equation*}
\mathbb  E \int^T_0 |Lu^n|^2 _{\Xi_2^\prime}\, dt \leq  2 \, \mathbb  E \int^T_0 | u^n|_2^2\, dt .
\end{equation*}
Multiplying both sides by $\epsilon^2 $, we obtain with  (\ref{uni1})
\begin{equation}\label{j4}
\mathbb  E |J ^n _4 |^2_{ H^1 (0, T; \, \Xi_2^\prime)} \mbox{ is bounded independently of  $n$ and $\epsilon$}.
\end{equation}

For  $J^n_6$,  Lemma \ref{lemitofrac} implies that, $\forall \,\, \alpha < \dfrac {1}{2}$,
\begin{align*}
\mathbb   E \left|\int^t _0  \sigma^n(u^n(s) ) \, d W(s)\right|^6_{W^{ \alpha, 6} (0, T;\, L^2(\dom))} & \lesssim \mathbb  E \int^t _0  \left | \sigma^n(u^n(s) )\right |^6_{\L} ds\\
& {\lesssim} \mathbb  E \int^t _0  \left | \sigma(u^n(s) )\right |^6_{\L} ds\\
&\leq (\mbox{by (\ref{sigma})}) \\
&\leq c^\prime\,c_B ^ 6\, \mathbb E \int^t _0  (1 + |u^n|)^  6 ds.\\
\end{align*}
This together with   (\ref{uni1})  
implies that
\begin{equation}\label{j5b}
\mathbb E \left |J^n_6 \right  |^2_{W^{\alpha, 6} (0, T;\, L^2 (\dom ))} \mbox{ is bounded independently of $n$ and $\epsilon$}, \,\,\,\,\forall \,\, \alpha < \dfrac {1}{2}.
\end{equation}
Hence we obtain (\ref{nostochastictermmm}).
%
%
%
%
%
%
%
%
Collecting the  estimates (\ref{j2}) and   (\ref{j3})-(\ref{j5b}), we obtain
\begin{equation}\label{eu1}
 \mathbb E |u ^n  |_{W^{\alpha, 2} (0, T;\, \Xi_2^\prime)} \mbox{ is bounded independently of $n$ and $\epsilon$}, \quad \alpha < \dfrac {1}{2}.
\end{equation}
By (\ref{uni1b}) we deduce
\begin{equation}\label{asderied}
\mathbb E |u ^{ n} |_{L^2 (0, T; \, H^1_0 (\dom)) }  \mbox{ is bounded independently of $n$},
\end{equation}
but the bounds may depend on $\epsilon$. From (\ref{eu1}) and (\ref{asderied}) we obtain  (\ref{eub}).

Observing from  (\ref{eq10-111intb}) that $u^n (t) - \int^t _0  \sigma^n(u^n  ) \, d W(s) = J^n_1 +  J^n_2 +  J^n_3 +  J^n_4 +  J^n_5 $, and applying (\ref{j2}), (\ref{j3}) and (\ref{j4}), we obtain (\ref{nostochastictermb}) as desired. \qed
\begin{rem}\label{uni}
See Lemma \ref{uniformesti} below for a variant of the proof of Lemma \ref{uniformestib} leading to the analogue of  bounds in (\ref{eub})-(\ref{nostochastictermmm}) but  independent of $\epsilon$. Note however that the proof in Lemma
\ref{uniformesti} for $u^\epsilon$ can not be applied here to $u^{\epsilon, n}$, because  multiplication  by  {$\sqrt{1+x}$ does not commute with $P^n$}, which prevents us from deducing  for now the estimates derived from (\ref{eq10-1111}) below.
\end{rem}

\noindent\emph{{Proof of  Theorem  \ref{uepsilon}.}}
The rest of the proof of Theorem  \ref{uepsilon} is classical (see  e.g. \cite{FlandoliGatarek1} and {\cite{DebusscheGlattHoltzTemam1}).
 Applying Lemma   \ref{lemcompact1} (of the Appendix) and Chebychev's inequality to the estimates (\ref{eub})-(\ref{nostochastictermmm}),   we can use   the same technic as that for the  proof of Lemma 4.1 in \cite{DebusscheGlattHoltzTemam1} to derive the compactness and tightness properties of the sequences   $(u^{\epsilon, n}(t), W(t))$ in $n$ for fixed $\epsilon$. Then we apply  the  Skorokhod embedding theorem   to construct some  subsequence   $\{(u^{\epsilon, n_k}(t), W(t))\}$  that converges strongly as $n_k \rightarrow \infty$, upon shifting the underlying probability basis. Then we pass to the limit on the Galerkin truncation (\ref{10-11g}) as $n_k \rightarrow \infty$ ($\epsilon$ fixed).
Note that we do not need to worry about passing to the limit on the boundary conditions, because they are all well-defined (and conserved) thanks to  (\ref{uni4b}).
Thus, we have established the existence of martingale solutions to the regularized stochastic ZK equation  (\ref{10-11b})-(\ref{10-2}),   (\ref{eq3}) and    (\ref{eq107}) in the sense of Definition \ref{defweaksolu}.

As for the  pathwise  solutions,
%
%
we first prove the pathwise uniqueness of  martingale solutions,
and then by  the Gy\"ongy-Krylov  Theorem we obtain the global existence of pathwise solutions  in the sense of Definition \ref{defstrongsol}.
%
To conclude, we have completed the proof of Theorem \ref{uepsilon}.  \qed

We will develop these steps below in more details in the more complicated case when $\epsilon \rightarrow 0$.

\section{Passage to the limit as $\epsilon \rightarrow 0$ to study the stochastic ZK equation} 
\label{passageto}

We now aim to study the stochastic solutions to the ZK equation basically by passing to the limit as $\epsilon \rightarrow 0$ in  (\ref{10-11b}) and the boundary conditions  (\ref{eq3}),  (\ref{eq107}), (\ref{10-4}) and  (\ref{10-2}).

\vskip 2 mm


\noindent\emph{Definition of solutions of the ZK equation.}
The definition of the martingale and pathwise solutions for the ZK equation are essentially the same as that for the regularized equation, with the necessary changes in the assumptions, equations and the function spaces.

\begin{deff}\label{defweaksoluu}
(Global Martingale Solutions) Let  $\mu _ {u_0}$ be  the probability measure of $u_0$ given as in (\ref{martingaleinitial})  on $L^2(\dom)$ and assume that (\ref{sigma}), (\ref{unifsigma}) and (\ref{f}) hold.

 A global martingale solution to the stochastic ZK equation (\ref{eq111})-(\ref{eq107}) and   (\ref{eq33}) (in the Dirichlet case) is defined as in   Definition \ref{defweaksolu} as a pair $(\tilde {\mathcal {S}}, \tilde u)$, such that
\begin{equation}\label{def}
 \tilde u \in L^6 ( \tilde \Omega; \,  L^\infty ( 0, T;\,  L^2(\dom)  )) \cap L^2 (\tilde \Omega; \, L^2 ( 0, T; \, H_0^1(\dom)  )),
\end{equation}
\begin{equation}\label{weakcont}
 \tilde u (\cdot, \omega) \in \mathcal C ([0, T];\,  L^2_w (\dom)))\,\,\tilde {\mathbb P} - a.s.,
\end{equation}
and $\tilde u$ satisfying  almost surely
\begin{equation}\label{defofsol}
 \tilde u (t)+ \int^{t}_0 (\Delta { \tilde u}_x  + c  { \tilde u}_x + { \tilde u }  { \tilde u}_x ) d s=  \tilde  u(0) + \int^{t}_0  f\, ds + \int^{t}_0 \sigma(  \tilde u) \, d \tilde W(s);
\end{equation}
the equality in (\ref{defofsol}) is understood {in the sense of distributions on $\mathcal D ( \dom$)} for every $ 0\leq t\leq T$.

Moreover $ \tilde u$ vanishes on $\partial \dom$ (since $ \tilde u \in L^2 (\tilde \Omega; \, L^2( 0, T; \, H_0^1(\dom)  ))$) and $  \tilde  u_x \big| _{x=1}=0$. For the latter, we observe that according to Lemma \ref{trace} below, $ \tilde  u_x \big| _{x=1}=0$ makes sense in a suitable  space for any $\tilde u$ satisfying (\ref{def}) and (\ref{defofsol}).
\end{deff}

\begin{deff}\label{defstrongsolll}
(Global Pathwise Solutions; Uniqueness)

Let $ {\mathcal S}:= ( \Omega, { \mathcal F}, \{ {\mathcal {F}}_t \}_{t\geq 0},  {\mathbb P}, \{ { W}^k \}_{k \geq 1} )$  be a  fixed stochastic basis and suppose that $u_0$ is an $L^2 (\dom)$-valued random variable (relative to $\mathcal S$) satisfying (\ref{u0}). We suppose that $\sigma$ and $f$ satisfy (\ref{sigma}),  (\ref{unifsigma}),  (\ref{sigmapa}) and (\ref{fp}).

\vskip 1 mm

(i) A global pathwise solution $u$ of (\ref{eq111})-(\ref{eq107}) and   (\ref{eq33}) is defined as in Definition  \ref{defstrongsol} with  (\ref{uepsiw})-(\ref{defsolb}) replaced by (\ref{def})-(\ref{defofsol}). Also
 note that $u$ vanishes  on  $\partial \dom$ (because $u$ $\in$ $L^2 (\Omega; \, L^2 ( (0, T); \, H_0^1(\dom)  )))$ and $u_x \big| _{x=1}=0$ which makes sense for  the same reasons as  for the martingale solution.

 \vskip 1 mm

(ii) Global pathwise uniqueness is defined in the same way as in Definition \ref{defstrongsol}.


%
%
%

\end{deff}

The strategy is the same as that in the case of the regularized stochastic ZK equation in Section \ref{lacalpe}: we  first derive the global existence of   martingale solutions, then prove the pathwise uniqueness of martingale solutions and hence deduce  the existence of  global pathwise solutions.

\subsection{Martingale solutions  in dimensions $2$ and $3$ }
\label{subsection12}
All the subsequent proofs are valid for $d=1$ or $2$, except  for (\ref{uni0b1p}) and   (\ref{uni4t}),  and for  the uniqueness in Section \ref{localpau}, which are only valid for $d=1$ (space dimension two).
\begin{thm}\label{existence}
 When $d=1$ or $2$, suppose that $\mu_0$ satisfies (\ref{martingaleinitial}), that $\sigma$ and $f$  maintain (\ref{sigma}), (\ref{unifsigma}) and (\ref{f}). Then there exists a global martingale solution $(\tilde {\mathcal S}, \tilde u )$ of (\ref{eq111})-(\ref{eq107}) and   (\ref{eq33}) in the sense of Definition \ref{defweaksoluu}.
\end{thm}
%

 Furthermore, when $d=1$, and if additionally    $f$ and $\sigma$ satisfy (\ref{fp}) and (\ref{sigmapa}), then the martingale solution  is  pathwise unique (see Proposition \ref{prop1} below).

To prove Theorem \ref{existence},  similar to the case of  the regularized stochastic ZK equation,   we first  derive the estimates leading to  weak convergence, then
using the  Skorokhod embedding theorem we   upgrade the weak convergence into the strong convergence,  with the probability basis shifted. Special measures will be taken  to pass to the limit  in the boundary conditions.

\subsubsection{Estimates and developments  independent of $\epsilon$.}\label{secunif}

We begin the proof of Theorem \ref{existence} by  deriving  the  estimates on $u^\epsilon$ valid as $\epsilon \rightarrow 0$. We observe that we can prove the  estimates  in (\ref{uni2})-(\ref{uni11}) under the new assumptions  in Theorem  \ref{existence}.

\begin{lem}\label{uniformest1} With the assumptions of Theorem \ref{existence},  when $d=1, 2$, we have the following  estimates valid as $\epsilon \rightarrow 0$:
\begin{equation}\label{uni21}
u_x^{\epsilon} \big|_{x=0} \mbox{ remains bounded in } L^2 (\Omega; \, L^2 (0, T; \, L^2 \left(I_{x^\perp}\right))),
\end{equation}
\begin{equation}\label{uni31}
\sqrt \epsilon u^{\epsilon} \mbox{ remains bounded in } L^2 (\Omega; \, L^2 (0, T; \, \Xi_1)),
\end{equation}
\begin{equation}\label{uni0b1}
 u^{\epsilon}\mbox{ remains bounded in } L^6 (\Omega; \, L^\infty (0, T; \, L^2 (\dom))),
\end{equation}

If we additionally assume that $u_0$ and $f$ satisfy (\ref{u0}) and  (\ref{fp}), then we have
\begin{equation}\label{uni0b1p}
 u^{\epsilon}\mbox{ remains bounded in } L^7 (\Omega; \, L^\infty (0, T; \, L^2 (\dom))).
\end{equation}

\end{lem}

\noindent \textbf{Proof.}
The estimates follow from (\ref{uni2})-(\ref{uni1}) (or (\ref{uni11})) by passing to the lower limit first in $n$ and then in $\epsilon$ using the lower semicontinuity  of the norms; indeed e.g. to show (\ref{uni0b1}), with (\ref{uni1}) we obtain $|u^\epsilon|_{L^6 (\Omega; \, L^\infty (0, T; \, L^2 (\dom)))}$ $\leq$ $\liminf_n |u^{\epsilon, n}|_{L^6 (\Omega; \, L^\infty (0, T; \, L^2 (\dom)))}$$\leq$$\kappa^\prime _1  $, for a constant $\kappa^\prime _1 $  independent of $\epsilon$.
\qed

\begin{lem}\label{uniformest}
The assumptions are those of  Theorem \ref{existence} with $d=1$ or $ 2$. We have the following  estimates valid as $\epsilon \rightarrow 0$:
\begin{equation}\label{uni4}
u^\epsilon \mbox{ remains bounded in } L^2 (\Omega; \, L^2 (0, T; \, H^1_0 (\dom))).
\end{equation}

 If furthermore we suppose that $u_0$ and $f$ satisfy (\ref{u0}) and  (\ref{fp}), and   $d=1$,  then we have
\begin{equation}\label{uni4t}
u^\epsilon \mbox{ remains bounded in } L^{7/2} (\Omega; \, L^2 (0, T; \, H^1_0 (\dom))).
\end{equation}
\end{lem}

\begin{rem}\label{pathu}
We will use (\ref{uni0b1p}) and   (\ref{uni4t}) only when dealing with the pathwise  uniqueness (see the calculations leading to (\ref{buu3}) below).
\end{rem}

\noindent \textbf{Proof of Lemma \ref{uniformest}.} The proof does not follow promptly from the estimates on the $u^{\epsilon, n}$ as that of (\ref{uni21})-(\ref{uni0b1p}), but they are derived directly from the solutions $u^\epsilon$ of the regularized equations; this is in fact the reason for which we introduced this regularization. Note that the solutions  $u^\epsilon$  are sufficiently regular for the following calculations to be valid.

We start by multiplying (\ref{eq10-111}) with $\sqrt {1+x}$, to find
 \begin{equation}\label{eq10-1111}
d ( \sqrt {1+x} \, u^\epsilon)=\sqrt {1+x}\, \mathcal N  (u^ \epsilon) \,d t + \sqrt {1+x} \, \sigma(u^\epsilon ) \, d W(t),
\end{equation}
where again $\mathcal N ( u^\epsilon) := - Au^\epsilon - B(u^\epsilon) - \epsilon \, L u^\epsilon  +f^\epsilon$.
Applying the  It\=o formula to (\ref{eq10-1111}), we obtain
 \begin{equation}\label{ito2}
 \begin{split}
 d|\sqrt {1+x}\,  u^\epsilon|^2 = & \, 2   \left (\sqrt {1+x}\, u^\epsilon , \,\sqrt {1+x}\,  \mathcal N ( u^\epsilon ) \right)\,dt\\
&+ 2\left (\sqrt {1+x}\,  u^\epsilon , \,\sqrt {1+x}\, \sigma (u^\epsilon)\,dW(t)\right )\\
&+   || \sqrt {1+x}\,  \sigma(u^\epsilon ) || ^2 _{\L }\, dt.
\end{split}
 \end{equation}
 We drop the super index $\epsilon$ for the moment and with exactly the same calculations  as in  the deterministic case (see the proof of Theorem 3.1 in \cite{SautTemamChuntian}), performed a.s. and for a.e. $t$, we have:
 \begin{equation}\label{eq10-30}
 \begin{split}
2   &\left (\sqrt {1+x}\, u , \,\sqrt {1+x}\,  \mathcal N ( u ) \right ) \\
& \quad\quad\quad= -|\nabla u|^{2}-2 | u_x |^{2} - (1-2\epsilon)\left| u_x \big|_{x=0} \right|^2_{L^2(I_{x^\perp})}   \\
& \quad\quad \quad -   2 \epsilon   \left( |\sqrt {1+x}  \, u_{xx}|^2 +  |\sqrt {1+x}  \, u_{yy}|^2 + |\sqrt {1+x}  \, u_{zz}|^2 )\right)\\
&\quad\quad  \quad + 2(f, \,( {1+x})\,  u ) + \dfrac{2}{3}   \int_{\dom}u ^{3}  \,d\,\dom   + c |  u |^2.
\end{split}
 \end{equation}
 Integrating both sides of (\ref{ito2}) in $t$ from $0$ to $s$, $0 \leq s  \leq T$, we find with (\ref{eq10-30}) that when say $\epsilon \leq 1/4$,
\begin{equation}\displaystyle \label{eq21}
\begin{split}
 \int^s _0  |\nabla u^\epsilon |^{2}\, dt
&\leq    |\sqrt {1+x}\,  u^\epsilon _0|^2  + 2  \int_0^s   (f^\epsilon, \, (1+x)u^\epsilon)\, d t+   \dfrac{2}{3}  \int ^s _0  |u^\epsilon|^{3}_{L^3(\dom)}\, dt \\
& \,\,\,\,\,\,\,\,\,\,+  c  \int_0^s  |  u^\epsilon  |^2\, dt+ \int_0^s   || \sqrt {1+x} \, \sigma(u^\epsilon ) || ^2 _{\L} \, dt \\
&\,\,\,\,\,\,\,\,\,\,+ 2 \int_0^s    \left ((1+x) u^\epsilon , \, \sigma (u^\epsilon)\,dW(t) \right ).
\end{split}
\end{equation}
For the first term on the right-hand side, using $H^{1/2}(\dom) \subset L^3(\dom)$ in dimension three, we have
$
 |u^\epsilon|^3_{L^3(\dom)} 
\leq c^\prime |u ^\epsilon |^{3/2} |\nabla u ^\epsilon|^{3/2}
 \leq \dfrac{1}{4} |\nabla u ^\epsilon |^{2} + c^{\prime}|u^\epsilon |^6 $;
hence taking  expectations on both sides of (\ref{eq21}) and using H\"older's inequality, we obtain
\begin{equation}\label{eq21t}
\begin{split}
\dfrac{1}{2}\,  \mathbb E \int^s _0  |\nabla u^\epsilon |^{2}\, dt
\lesssim & \, 2\,  \mathbb E  |\,  u^\epsilon _0|^2  +  \mathbb E \int_0^s   |f^\epsilon|^2 \, d t +   c^{\prime} \, \mathbb E   \int_0^s |u^\epsilon |^6\, d t + c^\prime\\
&\,+ \mathbb E  \int_0^s   || \sqrt {1+x}\, \sigma(u^\epsilon ) || ^2 _{\L} \, dt.
\end{split}
\end{equation}
Here the stochastic term vanishes.
We find  with (\ref{uni0b1}) and (\ref{eq21t})
\begin{equation}\label{eq23}
\mathbb  E \int^T_0  |\nabla u^\epsilon |^{2}\, dt
\leq \kappa _5,
\end{equation}
for a constant $\kappa_5$ depending only on $u_0$, $f$, $T$ and $\sigma$, and   independent of $\epsilon$; this implies (\ref{uni4}).

Returning to (\ref{eq21}), when  $d=1$, we have
\begin{equation}\label{buuxu}
\begin{split}
\int_0^s   |u^\epsilon   |_{L^3(\dom)}^3 d t
&\leq ( H^{1/3}(\dom) \subset L^3(\dom)\,\,\mbox{in\,\,dimension\,\,2})\\
&\leq  \int_0^s   |u^\epsilon  |^{2} |\nabla u ^\epsilon | d t \,d\,\dom \\
&\leq c^\prime  \sup_{0\leq t \leq s} |u^\epsilon  (t)|^{4} + \dfrac{1}{3} \left (\int_0^s   |\nabla u ^\epsilon  | d t \right )^2.
\end{split}\end{equation}
Hence (\ref{eq21})  implies
\begin{equation}\label{eq24tp}
\begin{split}
\dfrac{1}{2} \int^s_0  |\nabla u^\epsilon |^{2}\,  dt
&\leq   2 | u^\epsilon _0|^2  +   \int_0^s   | f^\epsilon |^2 \, d t + c^{\prime} \int_0^s   | u^\epsilon |^2 \, d t +    c^{\prime}  \sup_{0\leq t \leq s} |u^\epsilon (t) |^{4}  + c^\prime\\
& \,\,\,\,\,\,\,\,\,\,\,+ 2 \int_0^s   \left ((1+x) u^\epsilon , \, \sigma (u^\epsilon)\,dW(t)\right ).
\end{split}\end{equation}
Taking the supremum over $[0, T]$,  raising  both sides to the power $7/4$, then taking expectations, we obtain with Minkowski's inequality and Fubini's Theorem:
\begin{equation}\label{eq24tpp}\begin{split}\displaystyle
 \dfrac{1}{2}\,  \mathbb E \left( \int^T_0  |\nabla u^\epsilon |^{2}\,  dt\right)^{7/4}
&\lesssim  \mathbb E | u^\epsilon _0|^{7/2}  +  \mathbb E \int_0^T  | f^\epsilon |^{7/2} \, d t \\
& \,\,\,\,\,\,\,+  \mathbb E \int_0^T   | u^\epsilon |^{7/2} \, d t +    \mathbb E \sup_{0\leq s \leq T} |u ^\epsilon (s)  |^{7}  + c^\prime\\
& \,\,\,\,\,\,\,+ 2 \, \mathbb E \left[   \sup_{0\leq s \leq T}\, \int_0^s  \left | \left ((1+x) u^\epsilon , \, \sigma (u^\epsilon)\,dW(t)\right ) \right |\right ]^{7/4}.
\end{split}\end{equation}
For  the stochastic term, we have
\begin{align*}
&\mathbb E \left[   \sup_{0\leq s \leq T}\, \int_0^s  \left | ((1+x) u^\epsilon , \, \sigma (u^\epsilon)\,dW(t))\right |\right]^{7/4}\\
&\quad\quad\quad\quad\quad  \leq \mathbb E    \sup_{0\leq s \leq T}\left | \int_0^s  ((1+x) u^\epsilon , \, \sigma (u^\epsilon)\,dW(t))\right |^{7/4}\\
&\quad\quad\quad\quad\quad \leq ( \mbox{by the   Burkholder-Davis-Gundy inequality (\ref{burkholder})} )\\
&\quad\quad\quad\quad\quad \leq c _1 \, \mathbb E\, \left [ \left (\int ^T_0 |u^\epsilon |^{2} \, ||\sigma (u^\epsilon)||^2 _{\L }  \, d  t \right ) ^ {7/8} \right ]\\
&\quad\quad\quad\quad\quad  \lesssim   \mathbb E \left [ \left (\int ^T_0 |u^\epsilon|^{4}\, d t \right ) ^ {7/8}\right ] + c^\prime.
\end{align*}
This together with (\ref{eq24tpp})  implies
\begin{equation}\label{eq24tppp} \begin{split}
 \dfrac{1}{2}\,  \mathbb E \left( \int^T_0  |\nabla u^\epsilon |^{2}\,  dt\right)^{7/4}
&\lesssim  \mathbb E | u^\epsilon _0|^{7/2}  +  \mathbb E \int_0^T  | f^\epsilon |^{7/2} \, d t +  \mathbb E \int_0^T   | u^\epsilon |^{7/2} \, d t\\
&\,\,\,\,\,\,\,\,+    \mathbb E \sup_{0\leq s \leq T} |u^\epsilon   (s)|^{7} +  \mathbb E  \int ^r_0 |u^\epsilon |^{4}\, d t  + c^{\prime }.
\end{split}\end{equation}
Hence  (\ref{eq24tppp}) and  (\ref{uni0b1p}) imply
\begin{equation}\label{eq23}
\mathbb E \left( \int^T_0  |\nabla u^\epsilon |^{2}\,  dt\right)^{7/4}
\leq \kappa _6,
\end{equation}
for a constant $\kappa_6$ depending only on $u_0$, $f$, $T$ and $\sigma$, and   independent of $\epsilon$; this implies (\ref{uni4t}). The proof of Lemma  \ref{uniformest} is complete. \qed

\vskip 2mm

\noindent\emph{Estimates  in fractional Sobolev spaces.}
%
\begin{lem}\label{uniformesti}
With the same assumptions as in Theorem \ref{existence} and $d=1, 2$, we have
\begin{equation}\label{eu}
 \mathbb E |u ^\epsilon  |^2_{\mathcal Y }  \leq \kappa_7,
\end{equation}
\begin{equation}\label{nostochasticterm}
\mathbb E \left|   u^\epsilon (t)  - \int^t _0  \sigma(u^\epsilon  ) \, d W(s) \right| ^2 _  {H^1  (0, T;\,\Xi_2^\prime)} \leq  \kappa _8,
\end{equation}
\begin{equation}\label{nostochastictermm}
\quad\,\,\,\,\,\,\,\,\mathbb E \left|\int^t _0  \sigma(u^\epsilon  ) \, d W(s)  \right|^2_{W^{\alpha, 6} (0, T;\, L^2 (\dom ))} \leq  \kappa _9, \,\,\,\,\,\,\forall \,\, \alpha < \dfrac {1}{2},
\end{equation}
where ${\mathcal Y }$ is defined as in (\ref{y_1}), and $\kappa_7$, $\kappa_8$ and $\kappa _9$ are independent of $\epsilon$.
\end{lem}
\noindent \textbf{Proof.} By repeating the proof of Lemma \ref{uniformestib} (see Remark \ref{uni}), we see that we can obtain for $u^\epsilon$ the estimates analog to (\ref{eub})-(\ref{nostochastictermmm})     independent of $\epsilon$. The only point is to derive the estimate of $\mathbb E |u ^{\epsilon} |_{L^2 (0, T; \, H^1_0 (\dom)) } $ being bounded independently of $\epsilon$ (see (\ref{asderied}) correspondingly). For that we just need the estimate (\ref{uni4}). Hence Lemma \ref{uniformesti} is proven.  \qed

%

\subsubsection{Compactness arguments for  $\{ (u^\epsilon, W)\}_{\epsilon >0}$ }

With these  estimates independent of $\epsilon$ in hand, we can   establish the compactness of the family $(u^\epsilon(t), W(t))$.
For this purpose we consider the following phase spaces:
\begin{equation}\label{phasespaceb}
\mathcal{ X}_{u} = L^2 (0, T; \, L^2 (\dom)) \cap \mathcal{C } (0, T;\, H^{-5} (\dom)), \quad \mathcal{ X}_W = \mathcal{C} (0, T; \,\mathfrak{U}_0), \quad \mathcal{X} =\mathcal{ X}_u \times \mathcal{ X}_W.
\end{equation}
We then define the probability laws of $u^\epsilon(t)$ and $W(t)$ respectively in the corresponding phase spaces:
\begin{equation}\label{couple1}
\mu^\epsilon_u(\cdot) = \mathbb{P}(u^\epsilon \in \cdot), 
\end{equation}
and
\begin{equation}\label{couple2}
\mu_W(\cdot) =\mu^\epsilon_W(\cdot)  =  \mathbb{P}(W \in \cdot).
\end{equation}
This defines a family of probability measures $\mu^\epsilon := \mu^\epsilon_u \times \mu^\epsilon_W $ on the  phase space $\mathcal{  X}$. We now show that this family is tight in $\epsilon$. More precisely:
\begin{lem}\label{weakcom1}
We suppose that $d=1, 2$, and the hypotheses of Theorem \ref{existence} hold.    Consider the measures $\mu^\epsilon $ on $ \mathcal {  X}$ defined according to (\ref{couple1}) and (\ref{couple2}). Then the family $\{\mu^\epsilon  \}_{\epsilon >0}$ is tight and therefore weakly compact over the phase space $\mathcal {  X} $.
\end{lem}
\noindent{\textbf{Proof}}.  We can use the same technic as in the proof of Lemma 4.1 in \cite{DebusscheGlattHoltzTemam1}.  The main idea is to apply Lemma \ref{lemcompact1} (of the Appendix) and  Chebychev's inequality to  (\ref{eu})-(\ref{nostochastictermm}).
\qed

\vskip 2 mm

\noindent \emph{Strong convergence as $\epsilon \rightarrow 0$.}
Since  the family of measures $\{ \mu^\epsilon \}$ associated with the  family $ (u^\epsilon(t), W(t))$ is weakly compact on $\mathcal{  X}$, we deduce that $ \mu^\epsilon $ converges weakly  to a probability measure $ \mu$ on $\mathcal{  X}$ up to a subsequence. We can  apply the Skorokhod embedding theorem (see Theorem 2.4 in \cite{ZabczykDaPrato1}, also \cite{Billingsley2} and \cite{Jakubowski}\footnote{particularly in \cite{Jakubowski}, the theorem applies to $\mathcal{  X}$ as a Polish space, that is, a separable completely metrizable topological space.})  to deduce  the strong convergence of a further subsequence, that is
:
\begin{prop}\label{exsofmart}
Suppose that $\mu_0$ is a probability measure on $L^2 (\dom)$ that satisfies {(\ref{martingaleinitial})}. Then there exists a probability space $(\tilde{\Omega}, \tilde{\mathcal{F}}, \tilde{\mathbb{P}}  ) $, and a subsequence $\epsilon_k$  of  random vectors $ (\tilde u^{\epsilon_k}, \tilde{W}^{\epsilon_k} )$ with values in  $\mathcal{ X}$ ($\mathcal X$ defined in (\ref{phasespaceb})) such that

(i) $(\tilde u^{\epsilon_k}, \tilde{W}^{\epsilon_k} )$ have the same probability distributions as $ (u^{\epsilon_k}, W^{\epsilon_k})$.

(ii) $(\tilde u^{\epsilon_k},  \tilde{W}^{\epsilon_k} )$ converges almost surely as $\epsilon_k \rightarrow 0$, in the topology of $\mathcal {  X}$, to an element $(\tilde u, \tilde W) \in \mathcal X$, i.e.
\begin{equation}\label{strongconv1}
\tilde u^{\epsilon_k} \rightarrow \tilde u \mbox { strongly  in }L^2 (0, T; \, L^2 (\dom))) \cap \mathcal{C } ([0, T];\, H^{-5} (\dom))\,\,a.s.,
\end{equation}
\begin{equation}\label{strongconv2}
 \tilde{W}^{\epsilon_k} \rightarrow  \tilde W \mbox { strongly  in }\mathcal{C} ([0, T]; \, \mathfrak{U}_0)\,\,a.s.,
\end{equation}
where $ (\tilde u, \tilde W)$ has the  probability distribution $\mu$.

(iii) $\tilde{W}^{\epsilon_k}$  is a cylindrical Wiener process, relative to the filtration $\tilde{\mathcal{F}}^{\epsilon_k}_t$, given by the completion of the $\sigma$-algebra generated by $\{(\tilde{u}^{\epsilon_k} (s), \tilde{W}^{\epsilon_k} (s)); \,\,s \leq t\}$.

(iv) For each fixed $\epsilon_k$,  $\tilde u^{\epsilon_k} \in   L^2 (\tilde \Omega; \, L^2 (0, T; \, \Xi_2))$. Moreover, all the statistical estimates on $ u^{\epsilon_k} $ are valid for $\tilde u^{\epsilon_k}  $, in particular,  (\ref{uni0b1}) and (\ref{uni4}) hold.

(v) Each pair $(\tilde{u}^{\epsilon_k} , \tilde{W}^{\epsilon_k} )$ satisfies (\ref{10-11b}) as an equation in $L^2 (\dom)$ a.s.,  and satisfies the boundary conditions (\ref{eq3}),  (\ref{eq107}), (\ref{10-4}) and (\ref{10-2}) thanks to (iv), that is, $\tilde{u}^{\epsilon_k} (t)$ is adapted to  $\tilde{\mathcal{F}}^{\epsilon_k}_t$, and
\begin{equation}\label{eq10-11111}
  \begin{cases}
&d \tilde{u}^{\epsilon_k}=(- A\tilde{u}^{\epsilon_k} - B(\tilde{u}^{\epsilon_k}) - \epsilon_k \, L \tilde{u}^{\epsilon_k}  +f^{\epsilon_k}   )\,d t + \sigma(\tilde{u}^{\epsilon_k} ) \, d \tilde{W}^{\epsilon_k}(t),\\
& \tilde{u}^{\epsilon_k} =0 \mbox{ on } \partial \dom,\,\,\,\,\,\,\,\,\,\tilde{u}^{\epsilon_k} _x\big|_{x=1}=0,\\
 &\tilde{u}^{\epsilon_k}  _{xx}\big|_{x=0} = \tilde{u}^{\epsilon_k}  _{yy}\big|_{ {y=\pm \frac {\pi}{2}}  }=\tilde{u}^{\epsilon_k}  _{zz}\big|_{ {z=\pm \frac {\pi}{2}}  }=0,\\
&\tilde{u}^{\epsilon_k}(0)= \tilde u^{\epsilon_k} _0.
\end{cases}
\end{equation}
%
%


\end{prop}


\noindent\textbf{Proof.}   (i) and (ii) follow directly from the  Skorokhod  embedding  theorem.


To prove (iv), we first observe that thanks to Lemma \ref{analyticset} (of the Appendix), the space ${ L^2 (0, T; \, \Xi_2)}$ is a Borel set in the space $\mathcal X_u$,
  and hence the  integration $\int _{ L^2 (0, T; \, \Xi_2)} | u |^ {2} \, d\, \mu^{\epsilon_k}_u (u)$   makes sense, and by (i) we have for each $\epsilon_k$,
\begin{align*}
\mathbb E | {u}^{ \epsilon_k}|^2 _{ L^2 (0, T; \, \Xi_2)}&= \int _{ L^2 (0, T; \, \Xi_2)} | u |^ {2} \, d\, \mu^{ \epsilon_k}_u (u)=
\tilde { \mathbb E} | \tilde {u}^{ \epsilon_k}|^2 _{ L^2 (0, T; \, \Xi_2)}
 < (\mbox{by (\ref{uepsiw})}) <\infty.
\end{align*}
In the same way we would prove that all estimates on $u^{\epsilon}$ are valid for $\tilde u^{\epsilon_k}$, particularly (\ref{uni0b1}) and (\ref{uni4}).

To prove (v),   we define
\begin{align*}
\tilde  M^ {\epsilon_k} :  =  \int^T_0  \left| \tilde   {u}^{ \epsilon_k} (t)   +       \int^t _0 A {\tilde u}^{\epsilon_k}+ B(\tilde u^{\epsilon_k})  +  \epsilon_k L  \tilde{u}^{ \epsilon_k}-  f ^ { \epsilon_k} ds
-  \tilde u^{\epsilon_k} (0) - \int^t _0   \sigma(\tilde {u}^{ \epsilon_k} )  \, d \tilde W^ {\epsilon_k}(s) \right |^2 \, dt;
\end{align*}
then we can use the exact same technique  in \cite{Bensoussan1} to prove $ \tilde {\mathbb E} \, \dfrac{ \tilde   M^ {\epsilon_k}}{1+ \tilde  M^ {\epsilon_k}}=0$. Hence we obtain  (\ref{eq10-11111}).
\qed

\subsubsection{Passage to the limit}
\label{passage}

Now equipped with the strong convergences in (\ref{strongconv1}), we can consider passing to the limit on the regularized equation (\ref{eq10-11111})$_1$ as $\epsilon_k \rightarrow 0$.   Note that (\ref{eq10-11111})$_1$ is the version of (\ref{10-11b}) provided by the Skorokhod  embedding  theorem.

Thanks to  (\ref{uni0b1}) and (\ref{uni4}), we deduce the existence of an element
  \begin{equation}\label{spacefortildeu}
   {\tilde u }\in  L^6 (\tilde \Omega; \, L^\infty (0, T; \, L^2 (\dom)))\cap  L^2 (\tilde \Omega; \, L^2 (0, T; \, H^1_0 (\dom))),
   \end{equation}
 and a subsequence still denoted as $\epsilon_k$     such that
\begin{equation}\label{weakconv11}
\tilde u^{\epsilon_k} \rightharpoonup  {\tilde u }\mbox{ weak-star in }L^6 (\tilde \Omega; \, L^\infty (0, T; \, L^2 (\dom))),
\end{equation}
 and
 \begin{equation}\label{weakconv12}
\tilde u^{\epsilon_k} \rightharpoonup{\tilde { u} }\mbox{ weakly in } L^2 (\tilde \Omega; \, L^2 (0, T; \, H^1_0 (\dom))).
\end{equation}
Fixing  $u^ \sharp \in \Xi_2 $, by (\ref{weakconv12}) and   (\ref{weakconv11}) we can pass to the limit in the linear terms.

For the nonlinear term, for every $u^\sharp \in \Xi_2 $, we write a.s. and for a.e. $t$:
%
%
%
%
%
\begin{equation}\label{buus}\begin{split}
\bigg| \int^t _0& \left (B(\tilde u ^{\epsilon_k}) -  B(\tilde u ),  u^\sharp  \right ) \, ds \bigg| \\
& = \dfrac{1}{2}\left|  \int^t _0     \left(  (\tilde u ^{\epsilon_k} -  \tilde u)(\tilde u ^{\epsilon_k} +  \tilde u) ,   u^\sharp_x  \right) \, ds   \right|\\
& \leq \dfrac{1}{2}  \int^t _0     \left|  \tilde u ^{\epsilon_k} -  \tilde u\right|  \left|   \tilde u ^{\epsilon_k}+  \tilde u \right| |   u^\sharp_x |_{L^\infty (\dom)}   \, ds   \\
&\leq (\mbox{by the same calculations  as in (\ref{linf})})\\
& \leq \dfrac{1}{2}  \int^t _0     \left|  \tilde u ^{\epsilon_k} -  \tilde u\right|   \left|   \tilde u ^{\epsilon_k} +  \tilde u \right|  |   u^\sharp |_{\Xi_2}   \, ds\\
&\leq \dfrac{1}{2}   |u^\sharp|_{\Xi_2}      \left( \int^T _0     |\tilde u ^{\epsilon_k} -  \tilde u| ^2 \, ds \right)^ {1/2}     \left ( \int^T _0     |\tilde u ^{\epsilon_k} +  \tilde u| ^2 \, ds\right )^ {1/2}   .
\end{split}\end{equation}
Thus with (\ref{strongconv1}) and (\ref{uni0b1}),  we deduce that
\begin{equation}\label{strongcovfbu}
 \int^t_0 \left (B(\tilde u ^{\epsilon_k})   , u^\sharp \right ) \, ds  \rightarrow   \int^t_0 \left ( B(\tilde u )   , u^\sharp \right )\, ds \quad \mbox{ for } a.e.\,\, (\tilde \omega, t ) \in \tilde \Omega\times(0, T).
\end{equation}
We next establish the convergence for the nonlinear term in the space $L^1 (\tilde \Omega \times (0, T))$. We calculate as in (\ref{bu1}),
\begin{equation*}\label{buuss}
\begin{split}
\mathbb   E \int^T_0 \left| \int^t _0 \left (B(\tilde u ^{\epsilon_k}) ,   u^\sharp \right ) \, ds \right|^{2}  dt
\lesssim \mathbb   E \int^T_0  | \tilde u ^{\epsilon_k}  |^{4}   |u^\sharp|^2_{\Xi_2}\, ds
   \lesssim |u^\sharp|^2_{\Xi_2} \mathbb  E \int^T_0 |\tilde u ^{\epsilon_k}  |^ 4 \, ds.
\end{split}\end{equation*}
Thus by (\ref{uni0b1}), we have
\begin{equation*}
\left \{  \int^t _0\left (B(\tilde u ^{\epsilon_k }) ,   u^\sharp \right) \, ds  \right  \}_{\epsilon _k>0} \mbox{ is uniformly integrable for all {$\epsilon_k$} in $L^1 (\tilde{\Omega} \times (0, T))$}.
\end{equation*}
Hence thanks to  the Vitali convergence theorem, we conclude that
\begin{equation}\label{strongconvbuu}
  \int^t _0 \left <B(\tilde u ^{\epsilon_k}) ,   u^\sharp \right > \, ds \rightarrow  \int^t _0 \left <B(\tilde u  ) ,   u^\sharp  \right> \, ds  \mbox{\,\,\, in $ L^1 ( \tilde \Omega) \times (0, T)$} .
\end{equation}


For the stochastic term,
%
by (\ref{strongconv1}) we obtain
\begin{equation}\label{strongconverbu}
|\tilde u^{{\epsilon_k}}  -  \tilde u   |^2 \rightarrow 0,\quad \mbox{for a.e.}\,\, (\tilde \omega , t) \in \tilde \Omega\times(0, T).
\end{equation}
Thus, along with
(\ref{unifsigma})   we deduce
\begin{equation*}
|\sigma  ( \tilde u ^{{\epsilon_k}} ) - \sigma  (\tilde u )|_{L_2 (\mathfrak{U}, H)} \rightarrow 0, \quad \mbox{for a.e.}\,\, (\tilde \omega , t) \in \tilde \Omega\times(0, T).
\end{equation*}
On the other hand, we observe that
\begin{align*}
 \sup _{\epsilon_k} \mathbb  E   \left(   \int^T_0  |\sigma   ( \tilde u ^{{\epsilon_k}} ) |^6_{L_2 (\mathfrak{U}, H)}\, ds \right)\lesssim   \sup _{{\epsilon_k}}\mathbb   E \left( \int^T_0 (1+ | \tilde u ^{{\epsilon_k}}    | ^6)\, ds   \right) ,
 \end{align*}
 where we  made use of
 (\ref{sigma}). We therefore infer from (\ref{uni1}) that $| \sigma ( \tilde u ^{{\epsilon_k}} )|_{L_2 (\mathfrak{U}, H)}$ is uniformly integrable for $\epsilon_k$ in $L^q (\tilde \Omega \times (0, T))$ for any $q \in [1, 6)$. With the Vitali convergence theorem we deduce that, for all such $q\in [1, 6)$,
 \begin{equation}\label{sigmaconv}
 \sigma ( \tilde u ^{{\epsilon_k}} ) \rightarrow \sigma (\tilde u) \mbox{ in } L^q (\tilde \Omega; \, L^q((0, T), L_2 (\mathfrak{U}, H))).
 \end{equation}
Particularly (\ref{sigmaconv}) implies the convergence in probability of $\sigma ( \tilde u ^{{\epsilon_k}} )$ in $L^2 ((0, T), L_2 (\mathfrak{U}, H))$. Thus, along with the assumption (\ref{strongconv2}), we apply Lemma \ref{brownianconv} (of the Appendix) and deduce that
\begin{equation}\label{sigmainpro}
\int^t _0 \sigma ( \tilde u ^{{\epsilon_k}})\, d \tilde{W} ^ {{\epsilon_k}} \rightarrow \int^t_0 \sigma(\tilde u ) \, d\tilde{W},  \mbox{  in probability in  $L^2 ((0, T); \, L^2 (\dom))$}.
\end{equation}
By the Vitali convergence theorem using the estimates involving (\ref{burkholder}) and  (\ref{sigmaconv}),  from (\ref{sigmainpro}) we infer a stronger convergence result:
\begin{equation}\label{sigmainpro1}
\int^t _0 \sigma ( \tilde u ^{{\epsilon_k}})\, d \tilde{W} ^ {\epsilon_k} \rightarrow \int^t_0 \sigma(\tilde u) \, d\tilde{W},  \mbox{  in } L^2 (\tilde \Omega;\, L^2 ((0, T); L^2 (\dom)).
\end{equation}
Hence we can pass to the limit in (\ref{10-11b}), and obtain (\ref{defofsol}) as an equation in   $\Xi_2^\prime$

{For the initial condition, since (\ref{strongconv1}) and (\ref{spacefortildeu})  imply that $\tilde u ^\epsilon  \in L^\infty (0, T;\, L^2 (\dom)) $ $\cap$ $\mathcal{C } ([0, T]; \, H^{-5} (\dom))$ a.s., hence $\tilde u ^\epsilon  $ is weakly continuous with values in $L^2 (\dom)$ a.s.;  then (\ref{weakcont}) follows.}

 Having shown that the limit $\tilde u $ almost surely satisfies  (\ref{defofsol}) in the sense of distributions on $\mathcal D ( \dom$), we want now to address the question of the boundary  conditions. We need to be more careful because of the lack of regularity (see Lemma \ref{trace} below).

\vskip 2 mm

\noindent \emph{Passage to the limit on the boundary conditions.}
Since $\tilde u \in L^2 (0, T; \, H^1_0 (\dom))$ a.s. (see (\ref{weakconv12})), we deduce that $\tilde u $ satisfies  the Dirichlet boundary conditions. Hence there remains to show that the boundary condition
\begin{equation}\label{eq3p}
\tilde u _x \big|_{x=1}=0,
\end{equation}
{is satisfied} almost surely. This boundary condition is the object of Lemma \ref{trace} below where we show that $ \tilde u _x \big|_{x=1}$ is well defined when $\tilde u  \in L^6 (\tilde \Omega; \, L^\infty ( 0, T; \, L^2(\dom)  ))$ $\cap$ $L^2 (\tilde \Omega; \, L^2 ( 0, T; \,H^1(\dom)  ))$,
  and satisfies an equation like  (\ref{defofsol}).
\begin{lem}\label{trace}
We assume that $\tilde u \in  L^6 (\tilde \Omega; \,  L^\infty ( 0, T; \, L^2(\dom)  )) \cap L^2 (\tilde \Omega; \, L^2 ( 0, T; \, H^1(\dom)  )) $ satisfies (\ref{defofsol}) almost surely in the sense of distributions on $\mathcal D ( \dom$), for every $ 0\leq t\leq T$. Then
\begin{equation}\label{eq10-18}
\tilde u_x, \tilde u_{xx} \in \mathcal C _x (I_x; \,  \mathcal B), \,\,\,\,\,\,\mbox{where } \mathcal B=   L^{5/4}(\tilde \Omega; \,\,   H^{-3}((0, T )\times I_{x^\perp}),
\end{equation}
and, in particular,
\begin{equation}\label{eq10-28}
\tilde u_x \big|_{x=0, 1} \mbox{ and } \tilde u_{xx} \big|_{x=0, 1},
\end{equation}
are well defined in $\mathcal B$.

\end{lem}
\noindent{\textbf{Proof}}. 
If $\tilde u$ almost surely satisfies   (\ref{defofsol}), then
 $\tilde{U}: = \int^t _0 \tilde{u}\, ds$ satisfies
\begin{equation}\label{eq111lt}
     \dfrac{\partial  \tilde{U}}{\partial t} + \Delta \dfrac{\partial \tilde{U}}{ \partial x  }+ c \dfrac{ \partial  \tilde{U} }{ \partial x  }  =F\,\,\,a.s.,
    \end{equation}
where $F:= \tilde u_0      -\int^t _0 B(\tilde{u})\, ds+ \int^t _0 f  \, ds  + \int^t _0  \sigma(\tilde{u} ) \, d \tilde{W}(s)$.\\

For the term $\int^t _0 B(\tilde{u})\, ds$,  we note that by  (4.10) in \cite{SautTemam},
\begin{equation*}\label{eq10-334}
  |\tilde{u} \tilde{u}_x|_{L^{9/8}(\dom)} \leq |\tilde{u} |^{2/3} |\nabla \tilde{u} |^{4/3}, \mbox{ for } a.e.\,\, t \mbox { and } a.s..
\end{equation*}
Hence we have a.s.
 \begin{equation}\label{buu5}\begin{split}
  \left |\int^t _0 B(\tilde{u})\, ds \right |^{5/4}_{ L^{5/4}(0, T; \, L^{9/8} (\dom))} &=  \int^T_0 |\tilde{u}\tilde{u}_x | ^ {5/4}_{L^{9/8}(\dom)} \, dt\\
 & \lesssim \int^T_0 \left(|\tilde{u} |^{5/6}  \right) ^ {6} + \left(  |\nabla \tilde{u} |^{5/3} \right) ^{6/5} dt\\
&\lesssim \int^T_0 |\tilde{u}|^ 5   +  |\nabla \tilde{u} |^2  \, dt.
 \end{split}\end{equation}
Since $\tilde u  \in L^6 (\tilde \Omega; \,  L^\infty ( 0, T;\, L^2(\dom)  )) \cap L^2 (\tilde \Omega; \, L^2 ( 0, T;\, H^1(\dom)  )) $, taking expectations on both sides of  (\ref{buu5}) we have
\begin{equation*}
\tilde {\mathbb  E} \left |\int^t _0 B(\tilde{u})\, ds \right |^{5/4}_{ L^{5/4}(0, T; \, L^{9/8} (\dom))} < \infty,
\end{equation*}
that is
\begin{equation}\label{buu55}
\begin{split}
&\int^t _0 B(\tilde{u})\, ds  \mbox{ belongs to }  L ^  {5/4} (  \tilde \Omega; \,   L^{5/4}(0, T; \, L^{9/8} (\dom))  ),\\
&\mbox{ \,and hence belongs to }    L^{5/4} (I_x; \, L^{5/4}( \tilde \Omega \times (0, T) \times  I_{x^\perp})).
\end{split}
\end{equation}

For the term $\int^t _0  \sigma(\tilde{u} ) \, d \tilde{W}(s)$, from  (\ref{sigmainpro1}) we deduce that
$ \int^t _0  \sigma(\tilde{u} ) \, d \tilde{W}(s) \in L^2 (\tilde \Omega;\,L^2 (0, T; \, L^2 (\dom))$.

Applying the above estimates, we obtain that
\begin{equation}\label{j314}
\begin{split}
&F  \mbox{ belongs to }  L ^  {5/4} (  \tilde \Omega; \,   L^{5/4}(0, T; \, L^{9/8} (\dom))  ),\\
&\mbox{\,\,\, and hence  belongs to  }   L^{5/4} (I_x; \, L^{5/4}( \tilde \Omega \times (0, T) \times  I_{x^\perp}))).
\end{split}\end{equation}
Hence Lemma \ref{lineartracet} (of the Appendix) applies with $p=5/4$ and $\mathcal E =  L^{5/4}( \tilde \Omega \times (0, T) \times  I_{x^\perp})$, and from (\ref{eq10-18t})  we have
\begin{equation}\label{forevery}
\tilde U_x\mbox { and }\tilde U_{xx} \mbox{ belong to } \mathcal C _x (I_x; \, L^{5/4} ( \tilde \Omega; \,  H^{-2}((0, T )\times (I_{x^\perp}))).
\end{equation}
 Since
$\tilde U_x(t) =\int^t _0 \tilde u_x\, ds$, we have $ \frac{d\,\tilde U_x(t) }{d\,t }=\tilde u_x (t)$; differentiation in time maps continuously  $ H^{-2}(0, T ) $ into $H^{-3}(0, T )$ and from (\ref{forevery}) we thus infer
(\ref{eq10-18}) and (\ref{eq10-28}).
%
\qed

\vskip 2 mm

We now need to show that the boundary condition $\tilde u^ {\epsilon_k} _x \big|_{x=1}=0$, ``passes to the limit" to imply (\ref{eq3p}).
The idea is to apply Lemma \ref{A.2t} (of the Appendix) to $\tilde{U}^{\epsilon_k}(t): = \int^t _0 \tilde{u}^{\epsilon_k} \, ds$.
Rewriting (\ref{eq10-11111}) in an integral form and rearranging, we obtain a.s.
\begin{equation}\label{eq10-1111int}
\begin{split}
\tilde{u}^{\epsilon_k} (t)   +       \int^t _0 \Delta   \tilde{u}^{\epsilon_k} _x  \, ds & +c \int^t _0 \tilde{u}^{\epsilon_k} _x  \, ds+    \epsilon_k \int^t _0 L  \tilde{u}^{\epsilon_k} \, ds\\
&=\tilde u ^ {\epsilon_k}_0    -\int^t _0 B(\tilde{u}^{\epsilon_k})\, ds+\int^t _0 f^{\epsilon_k}  \, ds  + \int^t _0  \sigma(\tilde{u}^{\epsilon_k} ) \, d \tilde{W}^{\epsilon_k}(s).
\end{split}\end{equation}
Hence for almost every $\tilde \omega$,   $\tilde{U}^{\epsilon_k}$ satisfies the linearized parabolic regularized equation:
\begin{equation}\label{traceconverge}\displaystyle
 \begin{cases}
 &\,\,\,\dfrac{\partial \tilde{U}^{\epsilon_k}}{\partial t}+ \Delta \dfrac{\partial \tilde{U}^{\epsilon_k}}{\partial x}+ c \dfrac{\partial \tilde{U}^{\epsilon_k}}{\partial x} +  \epsilon_k\, L \tilde{U}^{\epsilon_k}  =F^{\epsilon_k},\\
 &\,\,\,\tilde U^{\epsilon_k} \big|_{x=0}= \tilde U^{\epsilon_k} \big|_{x=1}= \tilde U_x^{\epsilon_k}\big|_{x=1}= \tilde U_{xx}^{\epsilon_k} \big|_{x=0}=0,
\end{cases}
\end{equation}
%
where $F^{\epsilon_k}: =  \tilde u ^ {\epsilon_k}_0  -\int^t _0 B(\tilde{u}^{\epsilon_k})\, ds
+\int^t _0 f^{\epsilon_k}  \, ds  + \int^t _0  \sigma(\tilde{u}^{\epsilon_k} ) \, d \tilde{W}^{\epsilon_k}(s).$

For the term $ \int^t _0 B(\tilde{u}^{\epsilon_k})\, ds  $, by the same calculations as those  leading to (\ref{buu5}), we infer from (\ref{uni0b1}) and (\ref{uni4}) that
\begin{equation}\label{j31}
\tilde {\mathbb  E} \left |\int^t _0 B(\tilde{u}^{\epsilon_k})\, ds \right |^{5/4}_{  L^{5/4}(0, T; \,L^{9/8} (\dom))}  \mbox{ is bounded independently of $\epsilon_k$}.
\end{equation}
By (\ref{sigmainpro1}) we deduce that  $ \int^t _0  \sigma(\tilde{u}^{\epsilon_k} ) \, d \tilde{W}^{\epsilon_k}(s)$ remains bounded in $L^2 (\tilde \Omega;\,L^2 ((0, T); \,L^2 (\dom))$.
Collecting all the previous estimates
we conclude that $\tilde {\mathbb  E} |F ^{\epsilon_k} |^{5/4}_{  L^{5/4}(0, T; \,L^{9/8} (\dom))}$   is bounded independently of $\epsilon_k$, and hence
\begin{equation*}
\begin{split}
F ^{\epsilon_k} \mbox{ is bounded  independently of $\epsilon_k$ in $L^{5/4} (I_x; \, L^{5/4}( \tilde \Omega \times (0, T) \times I_{x^\perp}))$}.
\end{split}
\end{equation*}
Applying  Lemma \ref{A.2t} (of the Appendix) with $p=5/4$, $\tilde {\mathcal E}=  L^{5/4}(\tilde \Omega \times (0, T)\times I_{x^\perp})$ and
 $\tilde {\mathcal B}=  L^2 (\tilde \Omega;\, H^{-1}_t (0, T;\, L^2(I_{x^\perp}))) + L^2 (\tilde \Omega;\, L^2_t (0, T;\, H^{-4}(I_{x^\perp}   )))+   L^{5/4}(\tilde \Omega \times (0, T)\times (I_{x^\perp}))$,
we deduce that  $\tilde U^{\epsilon_k} _x \big|_{x=1} $ converges to $\tilde U _x  \big|_{x=1} $ weakly in $ \tilde {\mathcal B}$. Hence
\begin{equation}\label{uxboundary}
 \tilde U _x \big|_{x=1} (t)
=0,
\end{equation}
a.s. and  for a.e. $t \in (0, T)$. Since $\tilde U _x \big|_{x=1} (t)
=\int^t _0 \tilde{u}_x  \big|_{x=1} \, ds$,
  thanks to the Lebesgue differentiation theorem, we infer from (\ref{uxboundary}) that
$\tilde{u}_x  \big|_{x=1}(t) =0$ a.s. and for a.e. $t \in (0, T)$.
Thus we have finished the proof of Theorem \ref{existence}. \qed

\subsection{Pathwise solutions in dimension $2$ ($d=1$)}
\label{lacalp}
 We aim to establish the  existence of  pathwise  solutions  when $d=1$, that is:

%
%
%

\begin{thm}\label{existenceofpathwise}
 When $d=1$, assume that, relative to a fixed stochastic basis $\mathcal S$, $u_0$ satisfies (\ref{u0}),   and that $\sigma$ and $f$ satisfy (\ref{sigma}), (\ref{unifsigma}), (\ref{sigmapa}) and (\ref{fp}) respectively.   Then there exists a unique global pathwise solution $u$ which satisfies  (\ref{eq111})-(\ref{eq107}) and   (\ref{eq33}) in the sense of Definition \ref{defstrongsolll}.
\end{thm}

To prove this theorem,  we first establish the pathwise uniqueness of  martingale solutions and then apply the Gy\"ongy-Krylov  Theorem (Theorem \ref{krilov} of the Appendix). The difficulty lies in  deducing the pathwise uniqueness  due to a lack of regularity of the martingale solutions (see (\ref{def}) and (\ref{weakcont})). Adapting the idea from the deterministic case (see \cite{SautTemamChuntian}), we introduce a preliminary result concerning the existence and uniqueness of global pathwise solutions to the  linearized stochastic ZK equation with additive noise. More importantly, we establish an  energy inequality, which leads to a suitable estimate  of the difference of the solutions     for the application of the version  of the stochastic Gronwall lemma given in Lemma \ref{lemsg1} below.


\subsubsection{Linearized stochastic ZK equation with additive noise ($d=1$)}
\label{pathwisesolutions}

\begin{prop}\label{linear1}
When $d=1$, let $\mathcal S$ be a fixed stochastic basis, that is
 \begin{equation*}
 {\mathcal S}:= ( \Omega, { \mathcal F}, \{ {\mathcal {F}}_t \}_{t\geq 0},  {\mathbb P}, \{ { W}^k \}_{k \geq 1} ).
 \end{equation*}
We consider the linearized stochastic ZK equation ($c=0$),
\begin{equation}\label{eq111l}
\begin{cases}
  & d  \Ri + \Delta  \Ri _x \, d t=g\, dt + h \, d W(t),\\
   & \Ri(0)=\Ri  _0,
   \end{cases}
  \end{equation}
with the boundary conditions (\ref{eq3}) and (\ref{eq107}) for $\Ri$. We assume that
 \begin{equation}\label{rtau0}
  \Ri_0 \in L^2 (\Omega; \, L^2 (\dom)) ,
 \end{equation}
and $h$
and $g$ are given    predictable processes relative to the stochastic basis $\mathcal S$, such that
 \begin{equation}\label{eq4-3l}
g \in L^2 (\Omega; \, L^ {4/3}(0, T; \, L^{4/3}(\dom)))\cap L^2 (\Omega; \, L^2(0,T; \, \Xi_2^\prime)),
\end{equation}
and
\begin{equation}\label{eq4-31l}
h     \in L^2 (\Omega;   \,L^ 2(0, T; \, {\L })\cap L^2 (\Omega; \,  L^ 2(0, T; \, {L_2 (\mathfrak U, \, \Xi_1 )}).
\end{equation}
Then there exists a unique global pathwise solution $\Ri$ to (\ref{eq111l})
 which satisfies   (\ref{eq3}) and (\ref{eq107}), and
  such that
\begin{equation}\label{eq4-3lll}
 \Ri  \in L^2(\Omega; \,  L^ {\infty}(0, T; \, L^{2}(\dom)))    \cap L^2(\Omega;  \,L^ {2}(0, T; \, H^1_0(\dom)))  ,
\end{equation}
and
\begin{equation}\label{stronconl}
 \Ri (\cdot, \omega) \in \mathcal C ([0, T]; \,L^2_w (\dom)) \,\,a.s..
\end{equation}
Furthermore $\Ri$  satisfies the following  energy inequality
for any stopping time $\tau_b$ with $0 \leq \tau_b \leq T$,
\begin{equation}\label{eq16l2}
\begin{split}
 &\frac{1}{2} \,  \mathbb E\sup_{ 0 \leq s \leq \tau_b  } \, | \Ri (s) |^2 +  \mathbb  E \int^{\tau_b}_{0 } | \nabla  \Ri |^2 \, dt  \\
&\,\,\,\,\,\,\,\,\,\,\,\,\leq \mathbb E | \Ri(0)|^2 +  2 \, \mathbb E \int^{\tau_b}_{ 0} \left |(g,\,  (1+x)  \Ri  )\right|\, dt  + c ^{\prime   }\, \mathbb  E \int^{\tau_b}_{0 }\, || h || ^2 _{\L } \, dt.
\end{split}\end{equation}

\end{prop}

\noindent\textbf{Proof.}
%
We will first show the existence of the  solutions, which is similar to that of the nonlinear case, but only easier because the use of a compactness argument and the
derivation of strong convergence are not necessary for the linearized model. Then we will verify the uniqueness of the solutions, which  is direct  since the noise is additive. More precisely, the difference of two solutions satisfies a deterministic equation depending on the parameters $\omega \in \Omega$.
 Finally, we will deduce the energy inequality (\ref{eq16l2}) utilizing the duality between the spaces to which $g$ and $\Ri $ each belongs.

  We start by proving the existence of pathwise solutions with application of the parabolic regularization:
\begin{equation}\label{10-111l}\displaystyle
\begin{cases}
&d  \Ri ^\epsilon+ \left (\Delta  \Ri ^\epsilon _x+ \epsilon  L  \Ri ^\epsilon \right)dt =g^\epsilon\,dt + h\, d W(t),\\
&\Ri^\epsilon (0)=\Ri^\epsilon _0,
\end{cases}
\end{equation}
 supplemented with the  boundary conditions (\ref{eq3}),  (\ref{eq107}) and the additional boundary conditions (\ref{10-4}), (\ref{10-2}). As in Section \ref{regularize}, there exist  $\{\Ri^\epsilon _0\}_{\epsilon>0}$, a family of elements  in the space $L^2 (\Omega; \, \Xi_1) \cap  L^{22/3} (\Omega;\, L^2 (\dom))$ which are $\mathcal F_0$ measurable, and such that, as $\epsilon \rightarrow 0$,
\begin{equation}\label{initialcppp}
\Ri ^\epsilon_0 \rightarrow \Ri _0 \mbox{ in } L^2( \Omega;\,  L^2 (\dom)) \mbox{ strongly};
\end{equation}
and there exist
 $\{g^\epsilon\}_{\epsilon>0}  $,  a family of predictable processes relative to the stochastic basis $\mathcal S$, so that
 \begin{equation}\label{initialcpgg}
g^\epsilon \in  {L^\infty (\Omega; \, L^{22/3} (0, T; \, L^2 (\dom)))},
\end{equation}
\begin{equation}\label{initialcpg}
g ^\epsilon \rightarrow g \mbox{ in }  L^2(\Omega; \, L^ {4/3}(0, T; \, L^{4/3}(\dom))) \mbox{ strongly as } \epsilon \rightarrow 0.
\end{equation}

Since (\ref{eq4-31l}) corresponds to  (\ref{sigma})  and (\ref{unifsigma}), we can use a proof similar to that of
 Theorem \ref{uepsilon} 
 to deduce the existence and uniqueness of the global pathwise  solution $ \Ri ^\epsilon$ for each fixed $\epsilon$.  Note that for the proof of existence,  although $g^\epsilon$ depends on $\omega$, it will not be a problem for us; this is essentially because we can prove the existence of a $pathwise$ solution without referring to any compactness argument.

In the sequel, we will derive the estimates independent of $\epsilon$, then pass to the limit on the parabolic regularization, where  again we need to pay special attention to the boundary conditions.
 \vskip 2 mm

\noindent \emph{(i) Preliminary  estimates independent of $\epsilon$.}
We will prove the following  bounds on $\Ri^\epsilon $ as $\epsilon \rightarrow 0$:
\begin{equation}\label{uni1l}
 \Ri ^\epsilon   \mbox{ remains bounded in } L^2 (\Omega; \, L^\infty (0, T; \, L^2 (\dom))),
\end{equation}
\begin{equation}\label{uni4l}
 \Ri ^\epsilon    \mbox{ remains bounded in } L^2 (\Omega; \, L^2 (0, T; \, H^1_0 (\dom))).
\end{equation}

We start by multiplying both sides of (\ref{10-111l}) by $\sqrt {1+ x}$  and applying the It\=o formula, we find
 \begin{equation}\label{eqdul}\begin{split}
 d|\sqrt {1+ x}\, \Ri ^\epsilon|^2 &=2   \,(\sqrt {1+ x} \, \Ri ^\epsilon , \, \sqrt {1+ x}\,  Q (  \Ri ^\epsilon ) )\,dt \\
&\,\,\,\,\,\,\,+ 2(\sqrt {1+ x}\, \Ri ^\epsilon , \,\sqrt {1+ x}\, h\,dW(t))+   || \sqrt {1+ x}\, h || ^2 _{\L }\, dt,
\end{split} \end{equation}
 where $ Q(  \Ri ^\epsilon):= - \Delta \Ri_ x ^\epsilon  - \epsilon \, L  \Ri ^\epsilon  +g^\epsilon$.
Let some stopping times $\tau_a$, $\tau_b$ be given  so that  $0 \leq  \tau_a \leq \tau_b \leq T$; we integrate (\ref{eqdul}) from $\tau_a$ to $s$ and take the supremum  over $[ \tau_a ,  \tau_b  ]$. After taking expected values, and by  the same calculations as those leading to (\ref{eq21}),  we obtain that when $\epsilon \leq \dfrac{1}{4}$,
\begin{equation}\label{eq16l23}\begin{split}
\mathbb E\sup_{ \tau_a \leq s \leq {\tau_b}} \, | \Ri ^\epsilon (s) |^2 +
 \mathbb E \int^{\tau_b}_{\tau_a}   |\nabla  \Ri ^\epsilon |^{2}  dt&  \leq  2\,  \mathbb E\, | \Ri^\epsilon ( {\tau_a })|^2  +    2 \, \mathbb E \int^{\tau_b}_ {\tau_a }  \left|(g^\epsilon,\,   ({1+ x})\Ri ^\epsilon) \right|\,  dt\\
 &\,\,\,\,\, \,\,\,+ 2
\, \mathbb E \sup_{ {\tau_a }\leq s \leq {\tau_b}}     \int ^s_ {\tau_a } ((1+ x)\Ri^\epsilon,  h\,dW(t))\\
& \, \,\,\,\,\,\,\,+ 2\, \mathbb  E \int^{\tau_b}_ {\tau_a }  || h|| ^2 _{\L } \, dt.
\end{split}\end{equation}
For the  stochastic term, we have
\begin{align*}
\mathbb  E  \sup_{ {\tau_a }\leq s \leq {\tau_b}}  \left |\int ^s_ {\tau_a } ( (1+ x)\Ri^\epsilon , \, h\,dW(t)) \right|
&\lesssim (\mbox{by the   Burkholder-Davis-Gundy inequality} )\\
&\lesssim   \mathbb E\, \left [ \left (\int ^{\tau_b}_  {\tau_a } | \Ri ^\epsilon |^{2} \, ||h||^2 _{\L }  \, d  t \right ) ^ {1/2} \right ]\\
&\leq \dfrac{1}{4}\,  \mathbb  E   \sup_{{\tau_a}\leq s \leq {\tau_b}}| \Ri ^\epsilon |^{2} + c^{\prime } \, \mathbb  E  \int ^{\tau_b}_{\tau_a} || h|| ^2 _{\L }\, d  t.
\end{align*}
Hence (\ref{eq16l23}) implies
\begin{equation}\label{eq16lll}\begin{split}
\frac{1}{2}\, \mathbb E\sup_{ \tau_a \leq s \leq {\tau_b}}  \, | \Ri ^\epsilon (s) |^2  +
 \mathbb E \int^{\tau_b}_{\tau_a}   |\nabla  \Ri ^\epsilon |^{2}\,  dt &\leq 2\, \mathbb E\, | \Ri ^\epsilon( {\tau_a })|^2\\
 &\,\,\,\,\,\,\,\,+  2 \, \mathbb E \int^{\tau_b}_ {\tau_a } \left|  (g^\epsilon,\,   (1+x)\Ri ^\epsilon) \right|\, dt  \\
  &\,\,\,\,\,\,\,\,+c^\prime \,  \mathbb  E  \int ^{\tau_b}_{\tau_a} || h|| ^2 _{\L }\, d  t.
\end{split}\end{equation}
To estimate the term $\mathbb E \int^{\tau_b}_ {\tau_a } \left|  (g^\epsilon,\,   (1+x)\Ri ^\epsilon) \right|\, dt $, we
observe that a.s.
\begin{align*}
  \left| (g^\epsilon,\,   (1+x)\Ri ^\epsilon) \right| &\leq  |g^\epsilon|_{L^{4/3} (\dom)} |(1+x)\Ri ^\epsilon   |_{L^{4} (\dom)} \\
                             & \leq (\mbox{by Sobolev embedding in dimension $2$})\\
                             & \leq |g^\epsilon|_{L^{4/3} (\dom)} |\nabla \Ri ^\epsilon   |^{1/2}  | \Ri ^\epsilon   |^{1/2} \\
                              & \leq c^\prime |g^\epsilon|^{4/3}_{L^{4/3} (\dom)}  | \Ri ^\epsilon   |^{2/3} + \frac{1}{4} |\nabla \Ri ^\epsilon   |^{2} \\
                              & \leq c^\prime  |g^\epsilon|^{4/3}_{L^{4/3} (\dom)} \left(  | \Ri ^\epsilon   |^{2}  + 1  \right) + \frac{1}{4} |\nabla \Ri ^\epsilon   |^{2}.
\end{align*}
Applying the above estimates to (\ref{eq16lll}) we obtain
\begin{equation}\label{eq16lllll}\begin{split}
\frac{1}{2} \, \mathbb E\sup_{ \tau_a \leq s \leq {\tau_b}}\, | \Ri ^\epsilon (s) |^2 &+
 \frac{1}{2}  \,  \mathbb E \int^{\tau_b}_{\tau_a}   |\nabla  \Ri ^\epsilon |^{2}\,  dt\\
 &\,\,\,\,\leq 2\, \mathbb E\, | \Ri ^\epsilon( {\tau_a })|^2 +  2 c^\prime  \, \mathbb E \int^{\tau_b}_ {\tau_a }  |g^\epsilon|^{4/3}_{L^{4/3} (\dom)}  | \Ri ^\epsilon   |^{2}\, dt\\
 &\,\,\,\,\,\,\,\,\,\,\,+ \mathbb E \int^{\tau_b}_ {\tau_a }2 c^\prime |g^\epsilon|^{4/3}_{L^{4/3} (\dom)}  +c^{\prime }   || h|| ^2 _{\L }\, d  t.
\end{split}\end{equation}
Thanks to (\ref{initialcpgg}) and (\ref{eq4-31l}),   we can apply the stochastic Gronwall lemma (Lemma \ref{lemsg} below)
to (\ref{eq16lllll})
 to find
\begin{equation}\label{eq16llll}\begin{split}
\frac{1}{2}\,   \mathbb E\sup_{ 0 \leq s \leq T } \, | \Ri ^\epsilon (s) |^2 +
 \frac{1}{2}  \,  \mathbb E \int^{T}_{0} &  |\nabla  \Ri ^\epsilon |^{2}\,  dt\\
 & \lesssim \mathbb E\, | \Ri^\epsilon _0 |^2   + \mathbb E \int^{T}_ { 0 }  |g^\epsilon|^{4/3}_{L^{4/3} (\dom)}  +  || h|| ^2 _{\L }\, d  t.
\end{split}\end{equation}
Thanks to (\ref{initialcpg}), we have  $ |g^\epsilon| _{ L^2(\Omega; \, L^ {4/3}(0, T; \, L^{4/3}(\dom))) }   \leq |g| _{ L^2(\Omega; \, L^ {4/3}(0, T; \, L^{4/3}(\dom))) } + c^\prime$. Hence
$ \mathbb E \int^{T}_ { 0 }  |g^\epsilon|^{4/3}_{L^{4/3} (\dom)}  \lesssim    \mathbb E  \left[ \int ^ {T}  _0  |g^\epsilon| ^ {4/3} _{L^{4/3} (\dom)}\, dt\right ]^{3/2} + c^\prime \leq    |g| ^2 _{ L^2(\Omega; \, L^ {4/3}(0, T; \, L^{4/3}(\dom))) } + c^\prime $;  thus (\ref{eq16llll}) implies that
%
%
  \begin{align*}
\frac{1}{2} \,  \mathbb E\sup_{ 0 \leq s \leq T } \, | \Ri ^\epsilon (s) |^2 +
 \frac{1}{2}  \, \mathbb E \int^{T}_{0}   |\nabla  \Ri ^\epsilon |^{2}\,  dt\quad\quad \quad\,\,\,\,\,\quad\quad \quad\,\,\,\,\,\quad\quad \quad\,\,\,\,\,
\end{align*}
\begin{equation*}
\quad\quad \quad\,\,\,\,\, \quad\quad \quad\,\,\,\,\, \lesssim \mathbb E\, | \Ri _0 |^2   +  |g|^2  _{ L^2(\Omega; \, L^ {4/3}(0, T; \, L^{4/3}(\dom))) } + c^\prime +  || h|| ^2 _{\L }\, d  t.
\end{equation*}
    Hence  along  with (\ref{eq4-3l}) and    (\ref{eq4-31l}),     we obtain (\ref{uni1l}) and (\ref{uni4l}).

\vskip 2 mm

\noindent \emph{(ii) Estimates  in fractional Sobolev spaces.}
By the same proof as for  Lemma \ref{uniformesti}, with $\mathcal Y$ defined as in (\ref{y_1}), we  derive the following  estimates independent of $\epsilon$.
\begin{equation}\label{eul}
 \mathbb E |\Ri ^\epsilon |^2 _{\mathcal Y }  < \kappa_8,
\end{equation}
with $\kappa_8$ independent of $\epsilon$. This estimate will be useful to prove the continuity in time in (\ref{stronconl}).
%

\vskip 2 mm

\noindent \emph{(iii) Passage to the limit as $\epsilon_k \rightarrow 0$.}
With (\ref{uni1l}), (\ref{uni4l}) and (\ref{eul}), we deduce the following weak convergences, for a subsequence ${\epsilon_k} \rightarrow 0$:
\begin{equation}\label{weakconvl}
{ \Ri}^{\epsilon_k} \rightharpoonup  { \Ri} \mbox { weakly in }  L^2(   \Omega;  \,L^ {2}(0, T; \, H^1_0(\dom)))\cap L^2 (\Omega; \, W^{\alpha, 2} (0, T ;\, \Xi_2^\prime)) ,
 \end{equation}
 \begin{equation}\label{weakconvl1}
{ \Ri}^{\epsilon_k} \rightharpoonup    { \Ri} \mbox { weak star in } L^2(    \Omega; \,  L^ {\infty}(0, T; \, L^{2}(\dom))) .
 \end{equation}
We  can thus pass to the weak limit in (\ref{10-111l}) and obtain
\begin{equation}\label{defstrongsoll}
 \left <   \Ri (t ),  \Ri ^\sharp \right > +   \int^{t } _{ 0  }    \left <\Delta     \Ri _x - g  ,  \Ri^\sharp \right >\,ds    = \left  < \Ri_0,  \Ri ^\sharp  \right > +  \int^{t }_{0 } \left < h,  \Ri^\sharp \right  > \, d {W}  ,
\end{equation}
 for almost every $(  \omega, t) \in   \Omega \times (0, T)$ and every  $ \Ri ^ \sharp \in \Xi_2$.


To pass to the limit  on the boundary conditions (\ref{eq3}) and (\ref{eq107}), we use the same idea as in Section \ref{passage}. Firstly, we can prove an analogue of Lemma \ref{trace}; that is,   $ \Ri_x \big|_{x=1}$ is well defined if $ \Ri  \in L^2 ( \Omega; \, L^\infty ( 0, T; \, L^2(\dom)  ))$ $\cap$ $L^2 ( \Omega; \, L^2 ( 0, T; \, H^1(\dom)  ))$,
  and satisfies an equation like (\ref{defstrongsoll}). To show this, we just need to observe that thanks to    (\ref{eq4-3l}), Lemma \ref{lineartracet} applies with $p=4/3$ and $\mathcal E =   L^{4/3}(  \Omega \times (0, T) \times  I_{x^\perp})$. Secondly, we can pass to the limit on the boundary conditions applying  Lemma \ref{A.2t} (of the Appendix) with $p=4/3$, $\tilde {\mathcal E}=  L^{4/3}(  \Omega \times (0, T) \times  I_{x^\perp})$ and $ \tilde {\mathcal B}=  L^2 ( \Omega;\, H^{-1}_t (0, T;\, L^2(I_{x^\perp}))) + L^2 (\Omega;\, L^2_t (0, T;\, H^{-4}(I_{x^\perp}   )))+   L^{4/3}( \Omega \times (0, T)\times (I_{x^\perp}))$.

To prove (\ref{stronconl}), we infer from (\ref{weakconvl}) that $\Ri \in  W^{\alpha, 2} (0, T ;\, \Xi_2^\prime) \cap L^\infty (0, T; \, L^2 (\dom))\,\,a.s.$, and hence  $\Ri \in \mathcal C (0, T ;\, H^{-5} (\dom ))) \cap L^\infty (0, T; \, L^2 (\dom))$ a.s.. Thus $\Ri   $ is weakly continuous with values in $L^2 (\dom)$ almost surely, which implies    (\ref{stronconl}).

To conclude, we have proven  the existence of a global pathwise solution $   \Ri $  which satisfies (\ref{eq111l}),  (\ref{eq3}) and  (\ref{eq107}). 

\vskip 2 mm

\noindent \emph{(iv) Global pathwise uniqueness.}
     We assume that $\Ri_1, \Ri_2$ are two solutions of (\ref{eq111l}); setting $\Ri=\Ri_1-\Ri_2$,
we substract the equation (\ref{eq111l}) for $\Ri_1$ from that for $\Ri_2$;  we obtain that almost surely
\begin{align}\label{eq10-33333l}
\begin{cases}
\dfrac{\partial \Ri}{ \partial t} +  \Delta \Ri _x  =0,\\
\Ri_0=0.
\end{cases}
\end{align}
With  (\ref{eq4-3lll}),  we have $
 \Ri  \in  L^ {\infty}(0, T; \, L^{2}(\dom))   \cap  L^ {2}(0, T; \, H^1_0(\dom)) \,\,a.s.$. Hence   we can apply Lemma 3.2 in \cite{SautTemamChuntian} and  deduce that $\frac{d}{dt} | \Ri |^2 \leq 0$ for a.e. $\omega\in \Omega$ and $t \geq 0$; thus  $ \Ri (\omega)=0$ follows whenever $\Ri_0(\omega) =0$.



\vskip 2 mm

\noindent \emph{(v) Passage to the limit to obtain energy inequality (\ref{eq16l2}).}
From (\ref{eq16lll}), we obtain when $\tau_a =0$,
\begin{equation}\label{eq16llll2}\begin{split}
\frac{1}{2} \, \mathbb E\sup_{ 0 \leq s \leq {\tau_b}} \, | \Ri ^\epsilon (s) |^2  +
 \mathbb E \int^{\tau_b}_{0}   |\nabla  \Ri ^\epsilon |^{2}\,  dt& \leq 2 \, \mathbb E\, | \Ri^\epsilon _0 |^2+  2\,  \mathbb E \int^{\tau_b}_ {0 }   \left |(g^\epsilon,\,   (1+x)\Ri ^\epsilon)\right |\, dt \\
 &\,\,\,\,\,\,\,\,  +c^{\prime}\, \mathbb  E  \int ^{\tau_b}_{0} || h|| ^2 _{\L }\, d  t.
\end{split}\end{equation}
 We infer from (\ref{weakconvl}) and  (\ref{weakconvl1}) that for any $\tau_b$ with  $0 \leq \tau_b\leq T$,
\begin{equation*}\label{weakconvll}
{ \Ri}^{\epsilon_k} \textbf {1}_{  t \leq \tau_b } \rightharpoonup    { \Ri}  \textbf {1}_{ t \leq \tau_b }   \mbox { weakly in } L^2(   \Omega; \,  L^ {2}(0, T; \, H^1_0(\dom))),
 \end{equation*}
\begin{equation*}\label{weakconvll1}
{ \Ri}^{\epsilon_k} \textbf {1}_{  t \leq \tau_b } \rightharpoonup    { \Ri}  \textbf {1}_{ t \leq \tau_b }   \mbox {   weak star in  } L^2(    \Omega; \, L^ {\infty}(0, T; \, L^{2}(\dom))),
 \end{equation*}
and hence we can pass to the lower limit on the left-hand-side of (\ref{eq16llll2}).
{To pass to the limit on the term $ \mathbb E\, | \Ri^\epsilon _0 |^2 $, we use (\ref{initialcppp}). }

For the term $\mathbb E \int^{\tau_b}_ {0 } \left|  (g^\epsilon,\,   (1+x)\Ri ^\epsilon) \right| \, dt  $, we first note that in dimension $2$,
\begin{align*}
 \left( \int^T_{0} |\Ri^\epsilon |^4_{L^4 (\dom)} \, ds \right) ^{1/4} &\leq  \left( \int^T _{0} | \Ri^\epsilon |^2  |\nabla \Ri^\epsilon |^2  \, ds \right) ^{1/4} \\
& \leq \sup_{ 0 \leq s \leq T } \, | \Ri ^\epsilon  (s) |^{1/2}  \left ( \int^T _{0}  |\nabla \Ri^\epsilon |^2  \, ds  \right )^{1/4} \\
& \leq 2 \sup_{ 0 \leq s \leq T } \, | \Ri ^\epsilon  (s) | + 2 \left ( \int^T _{0}  |\nabla \Ri^\epsilon |^2  \, ds  \right )^{1/2}.
\end{align*}
Squaring both sides and taking the expectations we can use (\ref{uni1l}) and (\ref{uni4l})  to obtain that, as $\epsilon \rightarrow 0$,
\begin{equation}\label{uni1llll}
 \Ri ^\epsilon  \mbox{ remains bounded in } L^2 (\Omega; \, L^4 (0, T; \, L^4 (\dom))),
\end{equation}
and hence a subsequence of  $\Ri ^\epsilon $ converges weakly in the space $L^2 (\Omega; \, L^4 (0, T; \, L^4 (\dom)))$,
which is the dual of $L^2 (\Omega; \,L^{4/3} (0, T; \, L^{4/3} (\dom)))$.
 Since
 \begin{equation*}
  \mathbb E \int^{\tau_b}_ {0} \left|   (g^\epsilon ,\,   (1+x)\Ri ^\epsilon) \right|  dt  = \mathbb E \int^{T}_ {0} \left| \textbf{1}_{t \leq {\tau_b}} (g^\epsilon ,  \, (1+x)\Ri ^\epsilon)\right|  dt = \mathbb E \int^{T}_ {0} \left| (g^\epsilon \textbf{1} _{t \leq {\tau_b}} ,\,   (1+x)\Ri ^\epsilon)\right| dt,
   \end{equation*}we see that with  (\ref{initialcpg}),  $g^\epsilon\textbf{1} _{t \leq {\tau_b}} \rightarrow g\textbf{1} _{t \leq{\tau_b}} $ strongly in $L^2 (\Omega;\, L^{4/3} (0, T; \, L^{4/3} (\dom)))$, and hence the convergence of $ \mathbb E \int^{\tau_b}_ {0}  \left|  (g^\epsilon ,\,   (1+x)\Ri ^\epsilon)\right|\, dt  $      follows.

Thus we can pass to the lower limit on the left-hand side of  (\ref{eq16llll2})  and to the limit on the right-hand side of  (\ref{eq16llll2}),  and thus  deduce  (\ref{eq16l2}). Hence  we have completed the proof of Proposition \ref{linear1}.\qed

\subsubsection{Global pathwise uniqueness  for the full stochastic ZK equation ($d=1$)}
\label{localpau}
The following result establishes the pathwise uniqueness of martingale  solutions to (\ref{eq111})-(\ref{eq107}) and   (\ref{eq33}).
\begin{prop}\label{prop1}
When $d=1$, suppose that $(\tilde {\mathcal S}, \u)$ and $(\tilde {\mathcal S}, \v)$ are two global martingale solutions of (\ref{eq111})-(\ref{eq107}) and   (\ref{eq33}) , relative to the same stochastic basis.  We assume that the conditions imposed in Definition \ref{defstrongsolll} hold. We define
\begin{equation*}
\Omega _0 = \{ \u (0)  = \v (0)   \}.
\end{equation*}
Then $ \u$,  $\v$ are indistinguishable on $\Omega _0$ in the sense that
\begin{equation}\label{uniquness}
\tilde {\mathbb P} (  \textbf 1 _{ \Omega _0  } (  \u (t)  = \v (t) ) ) =1, \,\,\,\forall \,\,0\leq t \leq T.
\end{equation}
\end{prop}


\noindent\textbf{Proof.} We will mainly use (\ref{eq16l2}) from Proposition \ref{linear1} and the version of the stochastic Gronwall lemma given in Lemma \ref{lemsg1} below.
We define $ \Ri=\u - \v  $. Due to the bilinear term $B (\tilde u)$, when attempting to estimate ${\Ri}$, the terms that involve only $\u$ or $\v$ will arise. To deal with this issue we define the stopping times
\begin{equation}\label{stoppingtime}
\begin{split}
\tau^{(m)} & = \inf _{t\geq 0} \left \{ \sup_{s \in [0, t] }  |\u|^2 + \int ^t _0 | \nabla \u |^2 \, ds  + \sup_{s \in [0, t] }  |\v|^2 + \int ^t _0 | \nabla \v |^2 \, ds \geq  m \right \}\\
 &   =\sup _{t\geq 0} \left \{ \sup_{s \in [0, t] }  |\u|^2 + \int ^t _0 | \nabla \u |^2 \, ds  + \sup_{s \in [0, t] }  |\v|^2 + \int ^t _0 | \nabla \v |^2 \, ds \leq  m \right \}.
 \end{split}
\end{equation}
We  deduce from  (\ref{spacefortildeu}) that $\lim_{m\rightarrow \infty } \tau^ {(m) }= \infty  $.   Define $\Ra=  \textbf 1 _{ \Omega _0  } \Ri $, and  the  result will follow once we show that for any $ m $,
\begin{equation}\label{zero}
\tilde {\mathbb{E}} \left(   \sup_{[0, \tau^{(m)}\wedge T]}   | \Ra|^2   \right) =0.
\end{equation}

Subtracting the equation (\ref{eq111}) for $\v$ from that for $\u$, multiplying both sides by $ \textbf 1 _{ \Omega _0  }$, we arrive at the following equation for $\Ra $,
\begin{equation}\label{eq111llb}
\begin{cases}
& d   \Ra + \Delta   \Ra _x \, d t=\left(  - c  \Ra_x + \textbf 1 _{ \Omega _0  }( B(\v)  - B(\u) )     \right)\, dt + \textbf 1 _{ \Omega _0  }( \sigma (\u)  - \sigma( \v ) ) \, d \tilde  W(t),\\
 & \Ra(0)= 0.
\end{cases}
\end{equation}
Hence together with the stochastic basis $\tilde {\mathcal S}$,  we can regard  $\Ra $
  as a global pathwise solution      to (\ref{eq111llb}) written as (\ref{eq111l})  with   the boundary conditions  (\ref{eq3}) and (\ref{eq107}) , where   $g = - c  \Ra_x + \textbf 1 _{ \Omega _0  } ( B(\v)  - B(\u)) $ and $h =\textbf 1 _{ \Omega _0  } \left(   \sigma (\u)  - \sigma( \v  )\right)$.
%
To apply Proposition \ref{linear1}, now we only need to show that (\ref{eq4-3l}) and  (\ref{eq4-31l}) are satisfied. We  infer from (\ref{j33})
 that $g   \in L^2 (\tilde \Omega;  \, L^2(0,T; \, \Xi_2^\prime))$.
To show that $g$ also belongs to the space $L^2 (\tilde \Omega;  \, L^ {4/3}(0, T; \, L^{4/3}(\dom)))$,  we first note that
 ${ \Ri} _x  \in   L^2( \tilde \Omega;\,   L^2 ((0, T) ;\, L^2 (\dom) ))$.
 Next we estimate $B(u)$. By the Sobolev embedding theorem in dimension $2$,  we deduce that $|B(\u)|_{L^{4/3} (\dom)} \leq c^\prime|\u|_{L^4(\dom)}  |\u_x| \leq c^\prime   |\u|^{1/2}  | \nabla \u|^{1/2}|\u_x|$,  and hence almost surely
     \begin{equation}\label{buulinear}\begin{split}
 \left(\int ^ {T}  _0  |B(\u)| ^ {4/3} _{L^{4/3} (\dom)}\, dt\right )^{3/2} & \leq c _2 \left(  \int ^ {T}  _0   |\u|^{2/3}  | \nabla \u|^{2} \, dt \right )^{3/2}  \\
  & \leq c^\prime  \left(   \sup _{ t \in [0, T]} |\u|^{2/3} \,   \int ^ {T } _0   | \nabla \u|^{2} \, dt \right )^{3/2}  \\
     & = c^\prime   \sup _{ t \in [0, T ]} |\u| \left  ( \int ^ {T}  _0   | \nabla \u|^{2} \, dt \right )^{3/2}\\
&  \leq c^\prime   \sup _{ t \in [0, T ]} |\u|^7 + \left  ( \int ^ {T}  _0   | \nabla \u|^{2} \, dt \right )^{7/4}.
    \end{split}    \end{equation}
Since (\ref{uni0b1p}) and  (\ref{uni4t}) imply that  $\u$ and $\v$ both belong to the space $   L^{7} (\tilde \Omega; \,L^\infty ((0, T); \,H^1_0(\dom) ) )\cap    L^{7/2} (\tilde \Omega; \,L^2 ((0, T); \, H^1_0(\dom) ) ) $, taking expectations on both sides of (\ref{buulinear}) we obtain \begin{equation}\label{buu3}
       B(\u) \mbox{ belong to } L^2 (\tilde \Omega; \, L^{4/3}(0, T; \, L^{4/3}(\dom))).
\end{equation}
To conclude we infer  that $g \in L^2 (\tilde \Omega;  \, L^ {4/3}(0, T; \, L^{4/3}(\dom)))$.
 We infer from (\ref{sigmapa}) that
\begin{equation*}
|| h || _{L_2(\mathfrak U; \, \Xi_1) }  = || \sigma (\u) - \sigma (\v) ||  _{L_2(\mathfrak U; \, \Xi_1) }\leq   c_U   | \Ri|,
\end{equation*}
which implies that $h     \in  L^2 (\tilde\Omega; \,  L^ 2(0, T; \, {L_2 (\mathfrak U,  \Xi_1 )})$. Similarly, By  (\ref{unifsigma})
we   can deduce  that $h     \in L^2 (\tilde \Omega;   \, L^ 2(0, T; \, {\L }) $.
Thus we have proven that $h$ satisfies  (\ref{eq4-31l}).
%


Now Proposition \ref{linear1} applies, and we obtain
 (\ref{eq16l2})  for any $\tau_b$ with $0 \leq \tau_b \leq \tau^{(m)}\wedge T$, $\tau^{(m)}$ defined as in (\ref{stoppingtime}) (for notation simplicity, we will write $\tau^{(m)}:=\tau^{(m)}\wedge T$ from now on).
We then estimate the right-hand side of (\ref{eq16l2}). Thanks to  (\ref{stoppingtime}), we see that
 \begin{equation}\label{rxx}
  {\Ri (\cdot \wedge \tau^{(m)}  ) } \in  L^\infty ( \tilde \Omega; \,  L^\infty ((0, T) ;\, L^2 (\dom) )) \cap L^\infty (\tilde  \Omega; \,  L^2 ((0, T) ;\, H^1_0(\dom) )),
\end{equation}
 and hence the following calculations are all legitimate for $t \in (0, \tau^{(m)})$. We  observe that a.s. and for a.e. $t$:
\begin{equation}\label{buub}
\begin{split}
\left| (g,\,  (1+x)   \Ra ) \right|
&=  \left |- c( \Ra _x, \, (1+x )  \Ra   )+ (B(\v)  - B(\u)),   \, (1+x )   \Ra  )\right |\\
      &  = (\mbox{by } (\ref{buuuu}) \mbox{ of the Appendix})\\
      & = \left |\frac{c}{2}  | \Ra  |^2 +  \left (   {\Ra}^2, \, (1+x )\u_x  - \frac{1}{2} (  \v  +  (1+x) \v _ x)  \right ) \right|\\
 & \leq (\mbox{with } \gamma (t)  = |\u_x(t)|  +  |\v(t) |+  |\v_x(t)| )   \\
&\leq \frac{c}{ 2 }  |   \Ra   |^2  + c^{\prime}  \gamma(t) | \Ra |^2_{L^4(\dom)}  \\
& \leq (\mbox{by interpolation } H^{1/2} \subset L^4 \mbox{ in dimension }2)\\
&\leq \frac{c}{ 2 }  |\Ra   |^2  +   c^{\prime}  \gamma (t) |  \Ra | |\nabla { \Ra}  |\\
&\leq \frac{1}{2} |\nabla { \Ra}  |^2+ c^{\prime} \gamma^2(t) |  \Ra |^2 .
\end{split}
\end{equation}
Applying  (\ref{buub})  to  (\ref{eq16l2}), with  $\Ra( 0 )=0$ we obtain
\begin{equation}\label{eq16llllllll}
 \frac{1}{2}\,\tilde{  \mathbb E}\sup_{0\leq s \leq \tau_b  } \, |  \Ra (s) |^2 +  \dfrac{1}{2} \, \mathbb  E \int^{\tau_b}_{ 0 } |\nabla  \Ra |^2 \, dt  \leq  c^{\prime}\, \tilde{ \mathbb E} \int^{\tau_b}_{ 0 }\gamma^2(t)  |\Ra |^2\, dt,
\end{equation}
for any stopping time $\tau_b$ with $0 \leq  \tau_b \leq \tau ^ {(m)}$.
Along with (\ref{rxx}), the version of the stochastic Gronwall Lemma given in Lemma  \ref{lemsg1} below applies. Hence we  obtain
   (\ref{zero}).  This completes the proof of Proposition \ref{prop1}.
   \qed
%
%


Thanks to the pathwise uniqueness of  martingale solutions, we can  apply the Gy\"ongy-Krylov  Theorem to prove the existence of the pathwise solutions (for more details, see   \cite{DebusscheGlattHoltzTemam1}).

\vskip 2 mm

\noindent \textbf{Proof of Theorem \ref{existenceofpathwise}.}
We consider  the  families 
$(u^{\epsilon}, u^{\epsilon ^\prime}, W)$, where $u^{\epsilon}$ and  $u^{\epsilon ^\prime}$ are pathwise solutions to the parabolic regularization   (\ref{10-11b})-(\ref{10-2}),   (\ref{eq3}) and    (\ref{eq107}).
Then by (\ref{couple1}) and (\ref{couple2}), we can define the joint distributions of $(u^{\epsilon}, u^{\epsilon ^\prime}, W)$ as    $\nu ^ {\epsilon, \epsilon ^\prime} = \mu ^{\epsilon}_u \times \mu^{\epsilon^\prime} _u\times \mu _W$ on the phase space $\mathcal{X}_u \times \mathcal{X}_u \times \mathcal{X}_W$ ($\mathcal{X}_u $ and $\mathcal{X}_W$  defined in (\ref{phasespaceb})). With the same argument as for  Lemma \ref{weakcom1}, we can show that the family $\{ \nu ^ {\epsilon, \epsilon ^\prime} \}$ is tight  in $\epsilon$ and $\epsilon^\prime$.  By  the  Skorokhod  embedding  theorem,  we deduce the existence of a family      $(\tilde u^{\epsilon }, \tilde {\tilde u }^{\epsilon^\prime} , \tilde W)$ defined on a different probability space which converges almost surely to an element $ (\tilde u, \tilde {\tilde u }, \tilde W)$. By the same proof as for Proposition of \ref{exsofmart}, we can show that $(\tilde u^{\epsilon }, \tilde W)$ and $( \tilde {\tilde u }^{\epsilon^\prime} , \tilde W)$ both satisfy (i)-(v). In particular,     $(\tilde u^{\epsilon }, \tilde {\tilde u }^{\epsilon^\prime} )$  have the same probability distributions, $\mu ^{\epsilon}_u \times \mu^{\epsilon^\prime} _u$, as  $( u^{\epsilon }, {u }^{\epsilon^\prime} )$, and  the family  $\{ \mu ^{\epsilon}_u \times \mu^{\epsilon^\prime} _u\}_{\epsilon, \epsilon^\prime >0}$ is tight and hence converges weakly to a probability measure $\mu_1$, defined by $\mu_1  (\cdot ) = \mathbb P  (\tilde u, \tilde{\tilde u}\in \cdot )$. It is clear that $\tilde u$ and $\tilde{\tilde u}$ are both martingale solutions to   (\ref{eq111})-(\ref{eq107}) and   (\ref{eq33}), hence by the pathwise uniqueness (Proposition \ref{prop1}), $\tilde u=\tilde{\tilde u}$ a.s.. Thus
 \begin{equation*}
 \mu_1 (  \{ (u, v) \in  \mathcal{X}_u \times \mathcal{X}_u : u=v  \} )= \mathbb P (\tilde u=\tilde{\tilde u} \mbox{ in } \mathcal{X}_u   )= 1.
 \end{equation*}
We can apply
%
the Gy\"ongy-Krylov  Theorem  (Theorem \ref{krilov} of the Appendix) and deduce that the original family $u^{\epsilon}$ defined on the initial probability space converges in probability, and hence  converges
  almost surely up to a subsequence,
  to an element $u$  in the topology of $\mathcal{X}_u$. Thus we can pass to the limit on the regularized equation as explained in details in  Section \ref{passage}. To conclude we have  established the existence of a pathwise solutions to (\ref{eq111})-(\ref{eq107}) and   (\ref{eq33}), and we have completed the proof of Theorem \ref{existenceofpathwise}.  \qed

\begin{rem}
For the space periodic case, that is, (\ref{eq111}) and the boundary and initial conditions (\ref{eq3}),  (\ref{eq33}) and (\ref{eq1-26p}), the results will be the same with the Dirichlet case as discussed above. The reasoning will be similar as in \cite{SautTemamChuntian}.
\end{rem}

 \section{Appendix}
\label{sec:A}
In Section \ref{abstract} and  \ref{com}, we recall some results of deterministic nature. In Section \ref{sectrace} to \ref{gron} we present some results of probabilistic nature.

\subsection{Properties of $B(u)$}
\label{abstract}

In the article, we use the following properties of $ B(u)$ defined in (\ref{opera}).
\begin{lem}\label{lembuu}
Suppose $u, \,v\in H^1_0$,  then with  $\Ri = u-v$,  we have
\begin{equation}\label{buuuu}
 ((B(v)  - B(u),\,   (1+x) \Ri  ) =   (   {\Ri}^2, \, (1+x )u_x  - \frac{1}{2} (  v  +  (1+x) v _ x   ))  .
\end{equation}
\end{lem}
\noindent{\textbf{Proof}}.
\begin{align*}
 ((B(v)  - B(u),\,  (1+x)   \Ri ) & =   (\Ri   u _x, \, (1+x )  \Ri   ) +( v  \Ri_x,  \, (1+x )   \Ri  )\\
  &  =    ( \Ri^2, \, (1+x ) u_x ) + (  \Ri_x   \Ri, \,  (1+x) v   )\\
     &  =   (   {\Ri}^2, \, (1+x ) u_x ) - \frac{1}{2} ( { \Ri}^2 , \,  v  +  (1+x) v _ x   )\\
       & = (   {\Ri}^2, \, (1+x )u_x  - \frac{1}{2} (  v  +  (1+x) v _ x   )) .
      \end{align*}

\qed

\subsection{Compact embedding theorems}
\label{com}

We recall the theorems  from   \cite{FlandoliGatarek1} and \cite{Flandoli1} (see also \cite{Temam4}  for Lemma \ref{lemcompact1}).
\begin{deff}\label{fractional}
(The Fractional Derivative Space)
We assume that  $H$ is a separable Hilbert space. Given $\tilde p \geq 2$, $\alpha \in (0, 1)$, $ W^{ \alpha, \tilde p}(0, T;\, H )$ denotes the Sobolev space of all $h\in L^{\tilde p} (0, T; \, H)$ such that
\begin{equation}\label{deffract}
\int^T_0\int^T_0 \frac{ |  h(t) - h(s) |^ {\tilde p}_{H}  }{  |t-s| ^{1+\alpha \tilde p} }  \, dt\,ds < \infty,
\end{equation}
which is endowed with the Banach norm
\begin{equation}\label{fract}
 |h|  _  {W^{ \alpha, \tilde p}(0, T; \, H)  }  =  \left  (  \int^T_0 |h(t)|_H^{\tilde p} \, dt +       \int^T_0\int^T_0 \frac{ |  h(t) - h(s) |^ {\tilde p}_{H}  }{  |t-s| ^{1+\alpha \tilde p} }  \, dt\,ds  \right)^{1/ \tilde p} < \infty.
\end{equation}
\end{deff}

\begin{lem}\label{lemcompact1}
(i)
Let $\mathcal E_0 \subset \mathcal E\subset \mathcal {E} _1$ be Banach spaces, $\mathcal E_0 $ and $\mathcal {E} _1$ reflexive, with continuous injections and a compact embedding of $\mathcal E_0$ in $\mathcal E$. Let $1<p < \infty $ and $\alpha \in (0, 1)$ be given. Let $\mathcal Y$ be the space
\begin{equation}\label{space2}
\mathcal Y := L^p (0, T; \, \mathcal E_0 ) \cap  W^{\alpha, p} (0, T;\, \mathcal {E} _1),
\end{equation}
endowed with the natural norm. Then the embedding of $\mathcal Y$ in $L^p (0, T; \, \mathcal E)$ is compact.

\vskip 1 mm

(ii)
If $\mathcal {E}  \subset  \bar {\mathcal E}$ are two Banach spaces with $ \mathcal {E} $   compactly embedded in $ \bar {\mathcal E} $ ,  $1<p < \infty $ and $\alpha \in (0, 1)$  satisfy
\begin{equation*}
\alpha p > 1,
\end{equation*}
then the space  $W^{\alpha, p} (0, T;\, {\mathcal  E} )$ is compactly embedded into $\mathcal C ([0, T];\,  \bar {\mathcal E })$.
\end{lem}

 \subsection{Some generalized trace results}
\label{sectrace}
 The following trace result  is an extension of  the linear case of Lemma 3.1 of \cite{SautTemamChuntian}.
\begin{lem}\label{lineartracet}
  Let $u$ be a random process defined on a probability space $(\Omega, \mathcal F, \mathbb P)$. $F$ is a given function such that
\begin{equation}\label{spaceforg}
F\in  L^{p} (I_x ; \, \mathcal E ), \mbox{ where } \mathcal E = L^ {p} (  \Omega \times (0, T) \times I_{x^\perp}), \,\,p>1.
\end{equation}
   We assume that   $u \in L^ 2 (\Omega; \, L^2 (0, T; \, H^1 (\dom))) $
  satisfies almost surely  the following linear equation
\begin{equation}\label{new}
 u_ t +  \Delta   u_x+ c  u _ x    =F,
\end{equation}
that is, almost surely we have
\begin{equation*}
  u (t)+ \int^{t}_0 (\Delta {  u}_x  + c  {  u}_x ) d s=    u(0) + \int^{t}_0  F\, ds,
\end{equation*}
 {in the sense of distributions on $\mathcal D ( \dom$)} for every $ 0\leq t\leq T$ .
Then
\begin{equation}\label{eq10-18t}
u_x,  u_{xx} \in \mathcal C _x (I_x; \,\mathcal B),\,\,\, \mathcal B = L^2 (\Omega;\, H^{-2}( (0, T)\times I_{x^\perp}))\cap \mathcal E.
\end{equation}
and, in particular,
\begin{equation}\label{eq10-28t}
u_x \big|_{x=0, 1} \mbox{
and } u_{xx} \big|_{x=0, 1},
\end{equation}
are well defined in $  \mathcal B$.

\end{lem}
\noindent{\textbf{Proof}}. We write equation (\ref{new}) in the form
\begin{equation*}
u_{xxx}=F-cu_x- \Delta^\perp u_x -u_t.
\end{equation*}
Then clearly we have
\begin{equation}\label{eq10-5}
u_{xxx} \in L^{p \wedge 2} _x (I_x; \,L^{2 } ( \Omega; \, H^{-2}( (0, T)\times I_{x^\perp}  ) )\cap \mathcal E) ,\,\,\, p\geq 1,
\end{equation}
which implies that (\ref{eq10-18t}) holds.
\qed

We use Lemma \ref{lineartracet} in    the proof of Lemma \ref{trace} with $p=5/4$ and $\mathcal E =  L^{5/4}( \tilde \Omega \times (0, T) \times  I_{x^\perp})$, and
in the proof of Proposition  \ref{linear1} with $p=4/3$ and $\mathcal E = L^{4/3}(  \Omega \times (0, T) \times  I_{x^\perp})$.

Next we recall another trace result from \cite{SautTemamChuntian}.
\begin{lem}\label{A.2}
Let $\tilde {\mathcal B}$ be a reflexive Banach    space and let ${p\geq 1}$. We assume that two families of functions $u^\epsilon$ and $g^\epsilon \in L^p_x(I_x; \, \tilde {\mathcal B})  $ satisfy
\begin{equation*}
\begin{cases}
\quad u^\epsilon _{xxx} +\epsilon u^\epsilon_{xxxx}=g^\epsilon,\\
\quad u^\epsilon (0)= u^\epsilon (1)= u_x^\epsilon (1)= u_{xx}^\epsilon (0)=0,
\end{cases}
\end{equation*}
and that $g^\epsilon$ is bounded in $L^p_x(I_x;\,  \tilde {\mathcal B})  $  as $\epsilon \rightarrow 0$. Then $u^\epsilon _{xx}$ (and hence $u^\epsilon_x$, and $u^\epsilon$) is bounded in $L^\infty _x(I_x;\,  \tilde {\mathcal B})  $ as $\epsilon \rightarrow 0$. Furthermore, for any subsequences $u^ \epsilon \rightarrow u$ converging (strongly or weakly) in $L^q_x(I_x; \, \tilde {\mathcal B}) $, $  1<q $, $u^\epsilon_x (1)$ converges to $u_x (1)$ in $\tilde {\mathcal B}$ (weakly at least), and hence $u_x (1)=0$.
\end{lem}

We are now ready to prove the following trace result generalized  from the argument in \cite{SautTemamChuntian}.
\begin{lem}\label{A.2t}
 Let $\{u^\epsilon\}_{\epsilon >0}$ be a family of  random processes, all defined on the same probability space $(\Omega, \mathcal F, \mathbb P)$. We consider the following linearized  regularized equation
\begin{equation}\label{eq10-1112t}\displaystyle
 \begin{cases}
 &\,\,\, u^\epsilon_t+ \Delta  u^\epsilon _x+ c u^\epsilon _ x  + \epsilon L  u^\epsilon  =F^\epsilon,\\
 &\,\,\,u^\epsilon \big|_{x=0}= u^\epsilon \big|_{x=1}= u_x^\epsilon \big|_{x=1}= u_{xx}^\epsilon \big|_{x=0}=0,
\end{cases}
\end{equation}
where $F^\epsilon$ remains bounded in the reflexive Banach space $L^p (I_x,\tilde {\mathcal E})$, with  $\tilde {\mathcal E}:= L^p(\Omega \times (0, T)\times I_{x^\perp})$, $ p>1$. We suppose that  $u^\epsilon$ almost surely satisfies (\ref{eq10-1112t}),
 that is,  almost surely we have
\begin{equation*}
 u^\epsilon (t) +   \int^t _0    (\Delta  u^\epsilon_ x + c u ^\epsilon _ x  +  \epsilon L  u ^\epsilon- F^\epsilon ) \,ds    =  u^\epsilon (0),
\end{equation*}
in the sense of distributions on $\mathcal D ( \dom$) for every $ 0\leq t\leq T$ . We assume that $u ^\epsilon$ converges weakly to some $u$ in $L^2 (\Omega; \, L^2(0, T;\,H_0^1(\dom)))$ as $\epsilon \rightarrow 0$, then $u^\epsilon_x (1)$ converges   to $u_x (1)$ in $  {\tilde {\mathcal B}}$ specified  below,  and hence $u_x (1)=0$.
\end{lem}

\noindent{\textbf{Proof}}. By (\ref{eq10-1112t}) we have
 \begin{equation*}\label{traceconver}
u^\epsilon_{xxx} +\epsilon u ^\epsilon _{xxxx}  =F^\epsilon-u^\epsilon_t-cu^\epsilon_x-\Delta^\perp u^\epsilon _x - \epsilon u ^\epsilon _{yyyy}-\epsilon u^\epsilon _{zzzz}.
 \end{equation*}
We call the  right hand side $ g^\epsilon$. It is easy to observe that, since $u^\epsilon $ remains bounded in $L^2  (I_x; \, L^2 (\Omega \times (0, T) \times I_{x^\perp})$  as $\epsilon \rightarrow 0$,
then  $g^\epsilon $ remains bounded in the reflexive Banach space $ L_x^{p \wedge 2}  (I_x;\,  \tilde {\mathcal B}  ) $, $ p> 1$, where
\begin{equation}\label{mathcale}
\tilde {\mathcal B} =   L^2 (\Omega;\, H^{-1}_t (0, T;\, L^2(I_{x^\perp}))) + L^2 (\Omega;\, L^2_t (0, T;\, H^{-4}(I_{x^\perp}   )))+ \tilde { \mathcal E}.
\end{equation}
Thus we can apply Lemma \ref{A.2}  with this space $\tilde {\mathcal B}$ and obtain the convergence of the boundary term $u^\epsilon_x(1)$. \qed

Lemma \ref{A.2t} is applied in Section \ref{passage} with $p=5/4$ and  $\tilde {\mathcal E}=  L^{5/4}(\tilde \Omega \times (0, T)\times I_{x^\perp})$
and in the proof of Proposition  \ref{linear1},  with $p=4/3$, $\tilde {\mathcal E}=    L^{4/3}(  \Omega \times (0, T) \times  I_{x^\perp})$  and  $ \tilde {\mathcal B}=  L^2 ( \Omega;\, H^{-1}_t (0, T;\, L^2(I_{x^\perp}))) + L^2 (\Omega;\, L^2_t (0, T;\, H^{-4}(I_{x^\perp}   )))+   L^{4/3}( \Omega \times (0, T)\times (I_{x^\perp}))$.

\subsection{Convergence theorem for the noise term}
\label{converg}
The following convergence theorem for the stochastic integrals is used  to facilitate the passage to the limit in the parabolic regularization approximation. The statements  and proofs can be found in \cite{Bensoussan1},   \cite{GyongyKrylov1} and  \cite{DebusscheGlattHoltzTemam1}.
\begin{lem}\label{brownianconv}
Let $\{ \Omega, \mathcal F, \mathbb P \}$ be a fixed probability space, and $\mathcal X$ a separable Hilbert space. Consider a sequence of stochastic bases $\mathcal S_n:= ( \Omega, \mathcal F, \{\mathcal {F}^n_t \}_{t\geq 0}, \mathbb P,  W^n  )$, such  that each $W^n$ is a cylindrical Brownian motion (over $\mathfrak{U}$) with respect to $\{\mathcal {F}^n_t \}_{t\geq 0} $. We suppose that the $\{ {G^n} \}_{n\geq 1}$ are a sequence of $\mathcal X$-valued $\mathcal {F}^n_t $ predictable processes so that $G^n \in L^2 ((0, T);\,  L_2 (\mathfrak{U}, \mathcal X))$ a.s.. Finally consider $\mathcal S:= ( \Omega, \mathcal F, \{\mathcal {F}_t \}_{t\geq 0}, \mathbb P,  W )$ and a function $G$ $\in$ $L^2 ((0, T); \, L_2 (\mathfrak{U}, \,\mathcal X))$, which is $\mathcal {F}_t$ predictable. If
\begin{equation*}
W^n \rightarrow W \mbox{ in probability in } \mathcal C ([0, T]; \, \mathfrak U _0),
\end{equation*}
\begin{equation*}
G^n \rightarrow G \mbox{ in probability in }   L^2  ((0, T);\,  L_2 (\mathfrak{U}, \mathcal X)),
\end{equation*}
then
\begin{equation*}
\int^t _0 G^n \, d W^n  \rightarrow \int^t _0 G \, d W \mbox{ in probability in }  L^2  ((0, T); \, \mathcal X ) .
\end{equation*}
\end{lem}

\vskip 1 mm

Then we have the following lemma  based on the Burkholder-Davis-Gundy inequality and the notion of fractional derivatives in Definition \ref{fractional} (whose proof can be found in \cite{FlandoliGatarek1}).
\begin{lem}\label{lemitofrac}
Let $q \geq 2$, $\alpha > \dfrac{1}{2}$ be given so that $q \alpha >1$. Then for any predictable process $h \in L^{q} (\Omega \times (0, T); \, L_2(\mathfrak U, H))$, we have
\begin{equation*}
\int_0^t h(s) \, dW (s) \in L^{q } (\Omega; \, W^{ \alpha, q}(0, T;\,  H) ),
\end{equation*}
and there exists a constant $c^\prime=c^\prime  (q, \alpha) \geq 0$ independent of $h$ such that
\begin{equation}\label{lemitofraceq}
\mathbb E \left |\int^t _0  h(s) \, dW (s) \right |^{q}_{W^{ \alpha, q} (0, T;\, H)}\leq c^\prime(q,  \alpha)\, \mathbb E \int^t _0  | h(s) |^{q}_{L_2 (\mathfrak U, H)}\, ds.
\end{equation}
\end{lem}

\subsection{Some probability tools}

%
%

We  recall the   Gy\"ongy-Krylov  Theorem    from \cite{GyongyKrylov1},
which is used in proving the existence of pathwise solutions.
 \begin{thm}\label{krilov}
 Let $\mathcal X$ be a  {Polish} space. Suppose that $\{ Y_m \}$ is a sequence of $\mathcal X$-valued random variables on a probability space $(\Omega, \mathcal F, \mathbb P)$. Let $\{ \mu _{k,m } \}_{k, m \geq 1}$ be the sequence of joint laws of $\{Y_m  \}_{m \geq 1}$, that is
\begin{equation*}
\mu_{k,m }  (E):= \mathbb ((  Y_k, Y_m ) \in E), \quad E\in \mathcal B (\mathcal X\times \mathcal X).
\end{equation*}
Then
$\{ Y_m \}$ converges in probability if and only if for every subsequence of joint probability measures, $\{ \mu _{k_l, m_l}\}_{l \geq 0}$, there exists a further subsequence which converges weakly to a probability measure $\mu$ such that
\begin{equation}\label{diagno}
\mu (  \{  (u, v ) \in \mathcal X\times \mathcal  X : u=v  \}  )=1.
\end{equation}
\end{thm}

\subsection{The Jakubowski-Skorokhod representation theorem}

We  recall the following result from \cite{martinondrejat}.
\begin{lem}\label{analyticset1}
Let ${\mathcal A}_1 $ be a topological space such that there exists a sequence $\{f_m \}$ of continuous functions $f_m: {\mathcal A}_1 \rightarrow \mathcal R $ that separate points of ${\mathcal A}_1 $. Let ${\mathcal A}_2$ be a Polish space,  that is, {a separable} completely metrizable topological space,  and let $I:{\mathcal A}_2 \rightarrow  {\mathcal A}_1$ be a continuous injection. Then  $I(B)$ is a Borel set in ${\mathcal A}_1 $ whenever $B$ is Borel in ${\mathcal A}_2 $.
\end{lem}
The following result is a special case of Lemma \ref{analyticset1}.
\begin{lem}\label{analyticset}
Let  ${\mathcal A}_1 $ be a separable Hilbert space. Assume that ${\mathcal A}_2$ is a separable Hilbert space continuously injected into ${\mathcal A}_1 $.
 Then ${\mathcal A}_2$ is a  Borel set of ${\mathcal A}_1$.
\end{lem}
\noindent \textbf{Proof.} Firstly, it is clear that any separable Hilbert space ${\mathcal A}_1 $ satisfies the hypotheses of Lemma \ref{analyticset1}.
Since ${\mathcal A}_2$ is a  separable Hilbert space, hence it is a Polish space.
Now in Lemma \ref{analyticset1}, let $B$ be ${\mathcal A}_2$, which of course is a Borel set of ${\mathcal A}_2$. Then $I(B)=I ({\mathcal A}_2 ) =  {\mathcal A}_2$  is a Borel set in ${\mathcal A}_1 $   thanks to Lemma  \ref{analyticset1}. \qed

We  use Lemma \ref{analyticset}  in the proof of Proposition \ref{exsofmart}.

\subsection{The adapted stochastic Gronwall lemma}
\label{gron}

 We first recall the  stochastic Gronwall lemma from \cite{GlattHoltzZiane2} (see also \cite{MR}), then we present a variant result which is used in the proof of Proposition \ref{prop1}.
\begin{lem}\label{lemsg}
Fix $T > 0$. We assume that
\begin{equation*}
X, Y, Z, M: [0, T) \times \Omega \rightarrow \mathbb R,
\end{equation*}
are real valued, non-negative stochastic processes. Let $ 0 \leq \tau < T$ be a  stopping time so that
\begin{equation}\label{mx+z2}
\mathbb E \int ^{\tau} _{0} (M X + Z)\, ds < \infty.
\end{equation}
We suppose,  moreover that for some fixed constant $\kappa$,
\begin{equation}\label{kappa12}
 \int ^{\tau} _{0} M\, ds < \kappa, \quad a.s..
\end{equation}
Suppose that for all stopping times $\tau_a$, $\tau_b$ with  $ 0 \leq \tau_a \leq \tau_b \leq \tau$ we have
\begin{equation}\label{wehave2}
\mathbb E  \left(  \sup _{t \in [\tau_a, \tau_b]} X +    \int ^{\tau_b} _{ \tau_a } Y\, ds    \right)  \leq C_0  \,  \mathbb E \left( X ( \tau_a  )  + \int ^{\tau_b} _{ \tau_a }  (M X + Z)\, ds  \right),
\end{equation}
where $C_0$ is a constant independent of the choice of $\tau_a$, $\tau_b$. Then
\begin{equation}\label{res}
\mathbb E  \left(  \sup _{t \in [\tau_0, \tau]} X +    \int _{\tau_0} ^{ \tau} Y\, ds    \right)  \leq C  \,  \mathbb E \left( X (0  )  + \int ^{\tau} _{ \tau_0 }   Z \, ds  \right),
\end{equation}
where $C = C (C _0 , T, \kappa)$.
\end{lem}

When $X(0)=0$ and $Z=0$ we can weaken the hypotheses by requiring that (\ref{wehave2}) only holds for $\tau_a=0$ and all $\tau_b$, $ 0  \leq \tau_b \leq \tau$. We then obtain
\begin{lem}\label{lemsg1}
We assume that $X(0)=0$ and $Z=0$ in Lemma \ref{lemsg} and that (\ref{wehave2}) holds only for $\tau_a=0$ and all $\tau_b$, $ 0  \leq \tau_b \leq \tau$,  that is:
\begin{equation}\label{wehave23}
\mathbb E  \left(  \sup _{t \in [0, \tau_b]} X +    \int ^{\tau_b} _{ 0} Y\, ds    \right)  \leq C_0  \,  \mathbb E \left( X ( 0  )  + \int ^{\tau_b} _{ 0 }  M X \, ds  \right),
\end{equation}
where $C_0$ is a constant independent of the choice of $\tau_b$. Then the calculation (\ref{res}) holds true and reduces to
\begin{equation}\label{0}
\mathbb E    \sup _{t \in [0, \tau]} X = \mathbb E \, \int ^{\tau} _{ 0 } Y\, ds   =0.
\end{equation}
\end{lem}
\noindent \textbf{Proof.}
\textbf{Step 1.} We first show how to construct a finite sequence of stopping times
\begin{equation*}
0 < \tau_1 < ...< \tau_N < \tau_{N +1} =\tau\,\,\,\,a.s.,
\end{equation*}
so that
\begin{equation}\label{key}
\int^{ \tau_{j+1}} _{ \tau_{j} } M \, ds < \frac{1}{2C_0} \,\,\,a.s., \,\,\,\forall\,\, j =1, ..., N.
\end{equation}

We construct the sequence inductively. We start with time $0 $. We assume that $\tau_{j-1}$ is found. Then define
\begin{equation*}\label{stoppingtime22}
\tau_j :=\inf_{t\geq 0} \left \{  \int^{ t} _{\tau_{j-1} } M \, ds < \frac{1}{2C_0}  \right   \} \wedge \tau,
\end{equation*}
and $\tau_j >0$ is well-defined since $M > 0$ and it satisfies (\ref{kappa12}). Hence we have
\begin{equation}\label{suchthat2}
\int^{ \tau_j } _{\tau_{j-1} } M \, ds \geq  \frac{1}{2C_0}  , \quad \forall\, j \geq 1 \mbox{ such that } \tau_j < \tau.
\end{equation}
Now we claim that there exists a finite integer $N$ such that $\tau_N = \tau$, and
\begin{equation}\label{kappa22}
N \leq 2 C_0 \kappa + 1, \,\,\,a.s..
\end{equation}
We show this by contradiction; suppose that (\ref{kappa22}) is not true, then $N-1 > 2 C_0 \kappa$, and hence
\begin{align*}
\int^{ \tau_{N+1}} _{ 0 } M \, ds &=  \sum \limits_{j=1} ^ {N-1}  \int^{ \tau_{j+1} } _{\tau_{j} } M \, ds +   \int^{ \tau } _{\tau_{N} } M \, ds
\geq  \mbox{ with (\ref{suchthat2})} \geq (N-1 ) \frac{1}{2C_0} > \kappa.
\end{align*}
But this   contradicts with (\ref{kappa12}). Hence  (\ref{kappa22}) is proven, and we can choose the integer $N= \ulcorner  2 C_0 \kappa + 1   \urcorner$.

\textbf{Step 2.} Letting  $\tau_b =  \tau_1$ in (\ref{wehave23}), we have
\begin{equation}\label{wehave12}
\mathbb E  \left(  \sup _{t \in [0, \tau_1]} X +    \int ^{\tau_1} _{0} Y\, ds    \right)  \leq C_0 \, \mathbb E   \int ^{\tau_1} _{ 0 }  M X \, ds  ;
\end{equation}
from (\ref{wehave12}),   (\ref{mx+z2}) and (\ref{key}) we infer
\begin{equation}\label{eqx}
\mathbb E  \left( \dfrac{1}{2}  \sup _{t \in [0, \tau_1]} X +    \int ^{\tau_1} _{ 0} Y\, ds    \right)  =0.
\end{equation}
Thanks to (\ref{eqx}), for every   $\tau_b \geq \tau_1$ a.s., we find that
\begin{equation}\label{eqxxx}\displaystyle
\begin{split}
&\mathbb E   \sup _{t \in [0, \tau_b]} X = \mathbb E    \sup _{t \in [\tau_1, \tau_b]} X, \\
&\mathbb E \int ^{\tau_b} _{ 0 } Y\, ds   = \mathbb E\int ^{\tau_b} _{ \tau_1 } Y\, ds ,\\
 &\mathbb E\int ^{\tau_b} _{ 0 }  M X \, ds = \mathbb E\int ^{\tau_b} _{ \tau_1 }  M X \, ds.
 \end{split}
 \end{equation}
 Thus from (\ref{eqxxx}) and  (\ref{wehave23}), we infer that for every  $\tau_b \geq \tau_1$ a.s.,
\begin{equation}\label{wehave22m}
\mathbb E  \left(  \sup _{t \in [\tau_1, \tau_b]} X +    \int ^{\tau_b} _{ \tau_1 } Y\, ds    \right)  \leq C_0  \, \mathbb E\left( \int ^{\tau_b} _{ \tau_1 }  M X \, ds  \right).
\end{equation}
Setting $\tau_b=\tau_2$ in (\ref{wehave22m}), we have
\begin{equation}\label{eqxx}
\mathbb E  \left(  \sup _{t \in [\tau_1, \tau_2]} X +    \int ^{\tau_2} _{ \tau_1 } Y\, ds    \right)  \leq C_0  \, \mathbb E\left( \int ^{\tau_2} _{ \tau_1 }  M X \, ds  \right);
\end{equation}
with (\ref{eqxx}),  (\ref{mx+z2}) and (\ref{key}) we deduce
\begin{equation}\label{wehave1222}
\mathbb E  \left( \dfrac{1}{2}  \sup _{t \in [\tau_1, \tau_2]} X +    \int ^{\tau_2} _{ \tau_1} Y\, ds    \right)  =0.
\end{equation}
Hence by finite induction up to $N$ we obtain (\ref{0}).     \qed

\section*{Acknowledgements}
This work was partially supported by the National Science Foundation under the grants, DMS-0906440  and DMS 1206438, and by the Research Fund of Indiana University (Nathan Glatt-Holtz was partially supported by DMS - 1004638 and DMS-1313272).

  Both Nathan Glatt-Holtz and Chuntian Wang would like to acknowledge the Institute for Mathematics and its Applications at University of Minnesota where some of the work for this article was carried out.

Chuntian Wang  should like to thank   Professor Hari Bercovici for the discussions about the analytic sets, Professor Richard Bradley about the passage to the limit in the  terms involving stopping times,   Professor  Arnaud Debussche about the adaptivity of the solutions with the shifted probability basis due to the application of the    Skorokhod  embedding  theorem,  Professor Franco Flandoli for the explanation of some details of the proof in \cite{Bensoussan1} and the Jakubowski-Skorokhod representation theorem, and Professor Jean-Claude Saut for the suggestions in bibliographic  references.

   \newpage

\footnotesize
\bibliographystyle{amsalpha}

\bibliography{ref-3}

\normalsize

%

\end{document}